\documentclass[10pt,a4paper]{article}
\usepackage[all]{xy}
\usepackage[utf8]{inputenc}
\usepackage{amsmath,bm}
\usepackage{amssymb}
\usepackage{amsthm}
\usepackage{mathtools}
\usepackage{cases}
\usepackage{float}
\usepackage{here}
\usepackage{comment}
\usepackage{ulem}
\usepackage[T1]{fontenc}
\usepackage{textcomp}
\usepackage{graphicx, color}
\usepackage{tikz,epic,eso-pic}
\usepackage[hidelinks]{hyperref}

\theoremstyle{plain}
\newtheorem{thm}{Theorem}[section]
\newtheorem{lem}[thm]{Lemma}

\newtheorem{prp}[thm]{Proposition}

\theoremstyle{definition}
\newtheorem{df}[thm]{Definition}

\newtheorem{eg}[thm]{Example}

\newtheorem{rem}[thm]{Remark}

\numberwithin{equation}{section}
\allowdisplaybreaks[3]

\newcommand{\Z}{\mathbb{Z}}
\newcommand{\del}{\partial}

\title{
Magnitude, path and coarse homology
}

\author{
Yasuhiko Asao\\
\small Institute of Mathematics for Industry, Kyushu University\\
\small 744 Motooka, Nishi-ku, Fukuoka 819-0395, Japan\\
\small \texttt{asao@imi.kyushu-u.ac.jp}
}

\date{\empty}

\begin{document}
\maketitle

\begin{abstract}
In this article, we introduce {\it Borel--Moore} and {\it uniformly finite} versions of magnitude homology and path homology for locally finite graphs by allowing possibly infinite chains subject to local finiteness and boundedness conditions. In both settings, we construct magnitude--path spectral sequences extending the spectral sequence studied by the author \cite{A}. Their $E^1$-pages are given by the corresponding magnitude homologies, their $E^2$-pages contain the corresponding path homologies on the $p$-axis, and they converge to {\it coarse homology} and {\it uniformly finite homology}, respectively, in the sense of Nowak--Yu \cite{NY}. Moreover, the $E^{r+1}$-page is invariant under $r$-{\it homotopy}, whereas the homology to which the spectral sequence converges is invariant under coarse equivalence. Thus the spectral sequences interpolate between small-scale and large-scale geometry.

We also show that the Borel--Moore and uniformly finite magnitude homologies are products and restricted products, respectively, of ordinary magnitude homology groups with a fixed initial point. Hence the $E^1$-page is assembled from pointwise magnitude-homological data, providing a mechanism for combining local information into global large-scale invariants. As applications, we establish Mayer--Vietoris sequences for magnitude homology, path homology, the associated spectral-sequence pages, and coarse and uniformly finite homology. We also show that diagonality of ordinary magnitude homology implies degeneration at the $E^2$-page, yielding isomorphisms between path homology and the corresponding coarse homology. Finally, we discuss applications to Cayley graphs of graph products, including right-angled Artin groups and right-angled Coxeter groups.
\end{abstract}

\tableofcontents 
\section{Introduction}

One of the basic algebraic invariants in large-scale geometry  is {\it coarse homology} introduced by Roe \cite{Roe}, which we denote by $H_\ast$ in this paper. Its general construction involves some choices, namely a locally finite homology theory $h_\ast$ and a coarsening system $\{\mathcal{U}_i\}_i$ of the metric space. For a standard choice, coarse homology can be computed from {\it the Vietoris--Rips complexes} $V_R(X)$ of a metric space $X$ as $H_\ast(\varinjlim C_\ast^{\rm lf}(V_R(X)))$, where $C_\ast^{\rm lf}$ denotes {\it locally finite}, or {\it Borel--Moore} chain complex (see Remark \ref{BMhomology} for details). When $X$ is a locally finite graph, the  chain complex $\varinjlim C_\ast^{\rm lf}(V_R(X))$  admits a concrete description in terms of tuples of vertices :
\[
\varinjlim_{R\to\infty}
C_n^{\rm lf}(V_R(X))
=
\left\{ \sum_{\text{possibly\  infinite}} a_\sigma \sigma \mid \sigma = \{x_0, \dots, x_n\}, a_\sigma \in \Z, 
\sup_{a_\sigma\neq0}{\rm diam}\{x_0, \ldots, x_n\}<\infty
\right\}.
\]
 Up to the standard normalization and comparison between oriented simplices $\{x_0, \ldots, x_n\}$ and ordered tuples $(x_0, \ldots, x_n)$, this is the chain-level description of coarse homology used by Nowak--Yu \cite{NY} (see Remark \ref{Roedefinition} for details). Namely we can compute $H_\ast(X)$ as the homology of the chain complex $C_\ast(X)$ defined by 
\[
C_n(X)
=
\left\{ \sum_{\text{possibly\  infinite}}  a_\sigma \sigma \mid \sigma = (x_0, \dots, x_n), a_\sigma \in \Z, 
\sup_{a_\sigma\neq0}{\rm diam}\{x_0, \ldots, x_n\}<\infty
\right\},
\]
where we formulate this chain complex in terms of {\it functions on tuples} in this paper (see Section \ref{functionsontuples}).

Under {\it the bounded geometry} assumption (Definition \ref{bddgeom}), a bounded refinement of this theory is called {\it uniformly finite homology}, denoted by $H^{\rm uf}_\ast$ in this paper, introduced by Block--Weinberger \cite{blockweinberger}. In the tuple description above, its chain complex $C^{\rm uf}_\ast(X)$ is obtained by requiring, in addition, a uniform bound on the coefficients of the chains. It is remarkable that uniformly finite homology is closely related to amenability. In particular, Block--Weinberger showed that the vanishing of the zeroth uniformly finite homology characterizes non-amenability. 

Both of these coarse chain complexes $C_\ast(X)$ and $C_\ast^{\rm uf}(X)$, where the latter is defined when $X$ has bounded geometry, naturally carry filtrations reflecting the geometry of tuples. The most immediate one is the filtration by diameter, namely we can equip $C_n(X)$ with a filtration $F_\ell C_n(X)$ by replacing $<\infty$ with $\leq \ell$ in the definition of $C_n(X)$ above, and similarly for $C_\ast^{\rm uf}(X)$. In the present paper, however, we use instead the length of a tuple $(x_0,\ldots, x_n)$, defined by $L(x_0, \ldots, x_n)=\sum_{i=0}^{n-1} d(x_i, x_{i+1})$, and filter $C_\ast(X)$ by this quantity (Definition \ref{filteredchaincomplex}) as
\[
F_\ell C_n(X)
=
\left\{ \sum_{\text{possibly\  infinite}}  a_\sigma \sigma \mid \sigma = (x_0, \dots, x_n), a_\sigma \in \Z, 
\sup_{a_\sigma\neq0} L(x_0, \ldots, x_n)\leq \ell
\right\}.
\]

As usual in algebraic topology, a filtered chain complex gives rise to a spectral sequence. The main observation of this paper is that the spectral sequences arising from the length filtration have familiar homology theories on their early pages. In the setting of finite chains, the resulting spectral sequence, called {\it the magnitude--path spectral sequence}, was studied by the author in previous work \cite{A}, and subsequently investigated in detail by Hepworth--Roff \cite{HR}\cite{HRbigraded}, who introduced bigraded path homology as its full $E^2$-page and established several homological properties of the spectral sequence. More precisely, we introduce Borel--Moore and uniformly finite versions of {\it magnitude homology}, denoted by ${\rm MH}^{\rm BM}_{\ast, \ast}$ and ${\rm MH}^{\rm uf}_{\ast, \ast}$, and {\it path homology}, denoted by ${\rm PH}^{\rm BM}_{\ast}$ and ${\rm PH}^{\rm uf}_{\ast}$, for locally finite graphs (Definitions \ref{mag} and \ref{path}). The original versions of magnitude and path homology, ${\rm MH}_{\ast, \ast}$ and ${\rm PH}_{\ast}$, were introduced by Hepworth--Willerton \cite{HW} and Grigo'ryan--Lin--Muranov--Yau \cite{GLMY}, respectively.  Their chain complexes are defined using finite linear combinations of tuples of vertices (or allowed paths). Our Borel--Moore versions ${\rm MH}^{\rm BM}_{\ast, \ast}$ and ${\rm PH}^{\rm BM}_{\ast}$ are obtained by allowing possibly infinite linear combinations satisfying the same length or path conditions. When $X$ has bounded geometry, requiring in  addition a uniform bound on the coefficients gives the uniformly finite versions ${\rm MH}^{\rm uf}_{\ast, \ast}$ and ${\rm PH}^{\rm uf}_{\ast}$. Thus the homologies introduced here may be regarded as locally finite and uniformly finite extensions of the original magnitude and path homologies. With these definitions, we obtain the following, which extends \cite[Proposition~7.8]{A}.

\begin{thm}[Propositions \ref{weaklyconverge} and \ref{pathhomology}]
For a locally finite graph $X$, we denote the spectral sequences associated with the length filtrations on $C_\ast(X)$ and $C_\ast^{\rm uf}(X)$, where the latter is defined when $X$ has bounded geometry, by $E^{{\rm BM}, r}_{p, q}(X)$ and $E^{{\rm uf}, r}_{p, q}(X)$, respectively. Then they converge to $H_\ast(X)$ and $H_\ast^{\rm uf}(X)$, respectively, and we also have 
\begin{align*}
E^{{\rm BM}, 1}_{p, q}(X) &\cong {\rm MH}^{\rm BM}_{p+q, p}(X),\\
E^{{\rm uf}, 1}_{p, q}(X) &\cong {\rm MH}^{\rm uf}_{p+q, p}(X),
\end{align*}
for all $p, q$ with $p, p+q\geq0$, and
\begin{align*}
E^{{\rm BM}, 2}_{p, 0}(X) &\cong {\rm PH}^{\rm BM}_{p}(X),\\
E^{{\rm uf}, 2}_{p, 0}(X) &\cong {\rm PH}^{\rm uf}_{p}(X),
\end{align*}
for all $p\geq0$.
\end{thm}
In this way, magnitude homology and path homology appear as successive approximations to large-scale homological invariants. Moreover, as shown in Proposition \ref{invariance}, the $E^{r+1}$-page is invariant under $r$-{\it homotopy}, whereas the homology to which the spectral sequence converges is invariant under coarse equivalence. Thus the magnitude--path spectral sequence may be regarded as interpolating between small-scale and large-scale geometry. See also \cite{HRbigraded} for further homological properties of the pages of the magnitude--path spectral sequence.

Another important feature is that the Borel--Moore and uniformly finite magnitude homologies can be recovered directly from the ordinary magnitude homology with a fixed initial point. More precisely, let ${\rm MH}^x_{\ast,\ell}(X)$ denote the summand of the ordinary magnitude homology represented by tuples whose initial point is $x$. Then we have the following.

\begin{thm}[Proposition \ref{prod}]
For a locally finite graph $X$, we have 
\begin{align*}
{\rm MH}^{\rm BM}_{\ast, \ell}(X) &\cong \prod_{x \in X} {\rm MH}^{x}_{\ast, \ell}(X), \\
{\rm MH}^{\rm uf}_{\ast, \ell}(X) &\cong \overset{b}{\prod}_{x \in X} {\rm MH}^{x}_{\ast, \ell}(X),
\end{align*}
where the latter is defined when $X$ has bounded geometry, and $ \overset{b}{\prod}$ denotes the restricted product consisting of uniformly bounded families (Definition \ref{restrictedproduct}). 
\end{thm}
Thus the $E^1$-page is assembled from the magnitude homology groups ${\rm MH}^x_{\ast,\ell}(X)$ associated with individual initial points, whereas the spectral sequence eventually recovers global large-scale invariants. From this viewpoint, the spectral sequence provides a mechanism for assembling local magnitude-homological information into global geometric information. This product description also allows many structural properties of ordinary magnitude homology to be transferred directly to the Borel--Moore and uniformly finite counterparts. For example, diagonality of ordinary magnitude homology is equivalent to diagonality of Borel--Moore counterpart, and, under bounded geometry, of uniformly finite counterpart (Proposition \ref{diagonality}).
Similarly, using the description of the diagonal magnitude homology of graphs of girth at least $5$, we show that their Borel--Moore and uniformly finite path homologies vanish in all degrees at least $2$ (Proposition \ref{girth}).

The spectral sequences above provide useful tools for comparing coarse homology groups. Indeed, by the comparison theorem for convergent spectral sequences \cite[Theorem~5.2.12]{Wei}, a filtered chain map which induces an isomorphism on some finite page induces an isomorphism on the homology of the original chain complexes. Thus, in favorable situations, an isomorphism detected on a finite page of the spectral sequence can be promoted to an isomorphism of coarse or uniformly finite homology.

As a main application, we study Mayer--Vietoris sequences. Hepworth--Willerton \cite[Theorem~6.5, 6.6]{HW} proved a split Mayer--Vietoris sequence for magnitude homology for {\it a projecting decomposition} of a graph (see Section \ref{projectingdecomp}). We extend this result to the Borel--Moore and uniformly finite versions, and show that it induces a Mayer--Vietoris long exact sequence for path homology and, under suitable finiteness assumptions, split short exact sequences on the $E^r$-pages of the magnitude--path spectral sequences. We further prove a Mayer--Vietoris long exact sequence for coarse homology and uniformly finite homology. Namely we have the following.

\begin{thm}[Propositions \ref{BMMV}, \ref{PHMV}, \ref{EMV} and \ref{MVCO}]
Let $X$ be a locally finite graph, and let $(X; G, K)$ be a projecting decomposition. Then, we have the split exact sequence
\[
0 \to {\rm MH}^{\rm BM}_{n, \ell}(G\cap K) \xrightarrow{(\iota^G_{\ast}, -\iota^K_{\ast})}  {\rm MH}^{\rm BM}_{n, \ell}(G) \oplus  {\rm MH}^{\rm BM}_{n, \ell}(K) \xrightarrow{\iota_{G\ast}+\iota_{K\ast}}  {\rm MH}^{\rm BM}_{n, \ell}(X) \to 0,
\]
and long  exact sequences
\[
 \cdots \to {\rm PH}^{\rm BM}_{n+1}(X) \to {\rm PH}^{\rm BM}_{n}(G\cap K) \xrightarrow{(\iota^G_{\ast}, -\iota^K_{\ast})}  {\rm PH}^{\rm BM}_{n}(G) \oplus  {\rm PH}^{\rm BM}_{n}(K) \xrightarrow{\iota_{G\ast}+\iota_{K\ast}}  {\rm PH}^{\rm BM}_{n}(X) \to \cdots,
\]
and
\[
\cdots \to H_p(G\cap K) \xrightarrow{(\iota^G_{\ast}, -\iota^K_{\ast})}  H_{p}(G) \oplus  H_{p}(K) \xrightarrow{\iota_{G\ast}+\iota_{K\ast}}  H_{p}(X) \to H_{p-1}(G\cap K) \to \cdots.
\]
Moreover, if each fiber of the projecting map $\pi$ is finite, then we have the split short exact sequence
\[
0 \to E^{{\rm BM}, r}_{p, q}(G\cap K) \xrightarrow{(\iota^G_{\ast}, -\iota^K_{\ast})}  E^{{\rm BM}, r}_{p, q}(G) \oplus  E^{{\rm BM}, r}_{p, q}(K) \xrightarrow{\iota_{G\ast}+\iota_{K\ast}}  E^{{\rm BM}, r}_{p, q}(X) \to 0,
\]
for all $r\geq1$ and all $p, q$. Also, when $X$ has bounded geometry, the analogous statements hold for the uniformly finite versions. For the split short exact sequence on the $E^r$-pages, we further assume that the projecting map $\pi$ has uniformly finite fibers.
\end{thm}
The exact sequence for coarse homology is a special case of the general coarse Mayer--Vietoris principle for coarsely excisive decompositions \cite{Mit} (see Remark \ref{mitchener}). On the other hand, we do not know whether the corresponding Mayer--Vietoris sequence for uniformly finite homology has previously appeared explicitly in the literature.

Another useful situation occurs when the magnitude--path spectral sequences degenerate at the $E^2$-page. Recall that a graph $X$ is called {\it diagonal} if its ordinary magnitude homology is concentrated on the diagonal \cite[Definition~7.1]{HW}. As mentioned above, for locally finite graphs, this condition is equivalent to the diagonality of Borel--Moore magnitude homology and, under the bounded geometry assumption, also to the diagonality of uniformly finite magnitude homology. Consequently, for diagonal graphs the corresponding magnitude--path spectral sequences degenerate at the $E^2$-page, so that the Borel--Moore and uniformly finite path homologies compute coarse homology and uniformly finite homology, respectively. More precisely, we have the following.

\begin{thm}[Proposition \ref{e2deg}]
Let $X$ be a locally finite graph. If $X$ is diagonal in the sense of ordinary magnitude homology, then the magnitude--path spectral sequences $E^{{\rm BM}, \ast}_{\ast, \ast}$,  and, when $X$ has bounded geometry, also $E^{{\rm uf}, \ast}_{\ast, \ast}$, degenerate at the $E^2$-page, and we have
\begin{align*}
H_p(X)&\cong {\rm PH}^{\rm BM}_{p}(X), \\
H^{\rm uf}_p(X)&\cong {\rm PH}^{\rm uf}_{p}(X),
\end{align*}
for all $p\geq0$.
\end{thm}
This applies to every quasi-median graph (\S \ref{consequences}). Genevois \cite{Genevois} describes several classes of quasi-median graphs, including Cayley graphs of right-angled Artin groups and right-angled
Coxeter groups. Thus, when these graphs are locally finite, our spectral sequences degenerate at the $E^2$-page.

\paragraph{Organization of the paper.}
The remainder of this paper is organized as follows. In Section \ref{preliminary}, we collect preliminary material on normed abelian groups and restricted products, which is needed to handle uniformly finite chains, and recall the notion of a projecting decomposition of graphs. In Section \ref{functionsontuples}, we formulate infinite sums of tuples as functions on tuples, construct the corresponding simplicial objects, and introduce the length filtration. In Section \ref{MPC}, we define coarse and uniformly finite homology together with the Borel--Moore and uniformly finite versions of magnitude homology and path homology, and explain their relation to the classical definitions. In Section \ref{BMUFMPSS}, we construct the Borel--Moore and uniformly finite magnitude--path spectral sequences, study their convergence, functoriality and invariance properties, and give an example in which the spectral sequence does not stabilize at any finite page. We also discuss the application of $E^2$ degeneracy to coarse homology of quasi-median graphs. Finally, in Section \ref{MVS}, we establish the Mayer--Vietoris sequences discussed above for magnitude homology, path homology, the pages of the magnitude--path spectral sequences, and coarse and uniformly finite homology. 

\subsection*{Acknowledgements}
The author would like to thank Hitoshi Moriyoshi for pointing out the similarity between the coarse homology and the $E^\infty$-page of the magnitude--path spectral sequence. This observation was the starting point of this research. He also thanks Shun Wakatsuki for fruitful discussions. In particular, the proof of the Mayer--Vietoris sequence for $E^r$-pages is based on discussion with him. The author also thanks Tomoshige Yukita for introducing him to the theory of quasi-median graphs. The  author was supported by JSPS KAKENHI 24K16927.

\subsection*{Use of AI}
The author made extensive use of ChatGPT Pro (OpenAI) throughout the development of this work, both for mathematical discussions and for assistance with English writing. The mathematical discussions involved an iterative process in which the author proposed ideas and approaches, ChatGPT provided suggestions and possible arguments, and the author critically examined, modified, and further developed these suggestions through repeated exchanges. Some of the ideas used in the proofs emerged from this process, including suggestions originally made by ChatGPT. ChatGPT was also extensively used to improve the exposition, organization, and English expression of the manuscript.

All mathematical statements and proofs in the final manuscript have been checked and revised by the author, who takes full responsibility for the content.

\section{Preliminaries}\label{preliminary}
In this section, we collect some preliminary material used throughout the paper. We first introduce the category of normed abelian groups and study restricted products of such groups. In particular, we establish an $({\rm AB}4^\ast)$-type exactness property for restricted products under suitable uniform strictness assumptions. This will be used later to pass from magnitude homology with a fixed initial point to its uniformly finite version. We then fix our conventions for graphs and recall the notion of a projecting decomposition, which plays a central role in the Mayer--Vietoris arguments of Section~\ref{MVS}.
\subsection{Restricted products and an (AB$4^\ast$)-type property of normed abelian groups}\label{prelimabelian}
The abelian category ${\rm Ab}$ of abelian groups satisfies the following property, known as (AB$4^\ast$).
\begin{lem}\label{ab4}
For a family of short exact sequences of abelian groups
\[
0 \to A_i \xrightarrow{f_i} B_i \xrightarrow{g_i}C_i \to 0,
\]
the product $\prod_i$ preserves exactness, namely we have an exact sequence 
\[
0 \to \prod_iA_i \xrightarrow{ \prod_if_i}  \prod_iB_i \xrightarrow{ \prod_ig_i} \prod_iC_i\to 0.
\]
\end{lem}
In the following, we consider such a  (AB$4^\ast$)-type property in the category of {\it normed abelian groups}, where norms take values in non-negative integers.

\subsubsection{The preabelian category of normed abelian groups}
In this paper, {\it a normed abelian group} is an abelian group $A$ equipped with a map $|-| : A \to \Z_{\geq 0}$ satisfying 
\begin{itemize}
\item $|-a| = |a|, |a + b| \leq |a| + |b|$,
\item $|a| = 0$ if and only if $a = 0$.
\end{itemize}
If necessary, we denote the norm by $|-|_A$ indicating the base group. 
\begin{df}
Let $A$ and $B$ be normed abelian groups, and let $f : A \to B$ be a group homomorphism. 
\begin{enumerate}
\item We define
\[
|f| = \inf \{C>0 \mid \forall a, |f(a)|\leq C|a| \} \in [0, \infty],
\]
and we say $f$ is {\it Lipschitz} if $|f| < \infty$. Here and below, we use the convention that $\inf \emptyset = \infty$.
\item We define the category $|{\rm Ab}|_{\Z}$ whose objects are normed abelian groups and morphisms are  Lipschitz group homomorphisms.
\end{enumerate}
\end{df}

Note that isomorphisms in  the category $|{\rm Ab}|_{\Z}$ are exactly {\it bi-Lipschitz group isomorphisms}, namely a group isomorphism $f : A \to B$ such that there exist constants $c, C>0$ with $c|a|_A\leq |f(a)|_B\leq C|a|_A$ for all $a \in A$. The following is standard and we omit the proof.

\begin{lem}\label{lpnorm}
Let $(\Z^n, |-|_p)$ be a finite rank free abelian group with the $\ell_p$-norm for $p \in \{1, \infty\}$, namely we have $|(a_1, \dots, a_n)|_1= \sum_{k=1}^n|a_k|$ for $p =1$, and $|(a_1, \dots, a_n)|_\infty= \max_k |a_k|$. Then the identity on $\Z^n$ induces an isomorphism $(\Z^n, |-|_1) \cong (\Z^n, |-|_\infty)$ in $|{\rm Ab}|_{\Z}$. 
\end{lem}
In the following, we choose $\ell_1$ or $\ell_\infty$-norms for finite rank free abelian groups unless otherwise stated.

\begin{df}
For a subgroup $B \subset A$ of a normed abelian group $A$, we equip $B$ with a norm obtained by restricting that of $A$. Also, we equip the quotient group $A/B$ with the norm defined by
\[
|[a]|_{A/B} = \inf_{b\in B} |a+b|.
\]
This is indeed a norm; in particular, we have $|[a]| = 0$ if and only if $[a]=0$, since the norm takes values in $\Z_{\geq0}$.
\end{df}
Note that, for a subgroup $B\subset A$, the inclusion $B \to A$ and the quotient $A \to A/B$ are Lipschitz homomorphisms with respect to the norms defined above.

The category $|{\rm Ab}|_{\Z}$ is obviously pre-additive, namely it is ${\rm Ab}$-enriched. Moreover, we can show that this category is {\it pre-abelian}, namely it has biproducts, kernels and cokernels. In fact, we can define the biproduct of normed abelian groups $A\oplus B$ as the biproduct in ${\rm Ab}$ with the norm defined by $|(a, b)| = \max \{|a|, |b|\}$. Then the projections $A\oplus B \to A$ and $A\oplus B \to B$ are Lipschitz, and it is straightforward to verify the universalities of the product and the coproduct. Also, we can define ${\rm ker} f$ and ${\rm coker} f$ of a morphism $f : A\to B$ in $|{\rm Ab}|_{\Z}$ as groups ${\rm ker} f \subset A$ and $B/{\rm im}f$ with norms defined above.  The universality of the kernel is obvious, and that of cokernel is straightforward. 

\begin{lem}\label{zovera}
Let $A \subset \Z^m$ and $B\subset \Z^n$ be subgroups. Then every group homomorphism $f : \Z^m/A \to \Z^n/B$ is a morphism in $|{\rm Ab}|_{\Z}$. In particular, $f$ is an isomorphism in $|{\rm Ab}|_{\Z}$ if it is a group isomorphism.
\end{lem}
\begin{proof}
By Lemma \ref{lpnorm}, the $\ell^1$- and $\ell^\infty$-norms on $\Z^m$ are bi-Lipschitz equivalent. Taking the infimum over $A$, the corresponding quotient norms on $ \Z^m/A$ are also bi-Lipschitz equivalent, and similarly for $\Z^n/B$.  Hence we may choose the $\ell^1$-norm. In this case, the norms are identified with the word norms on the finitely generated groups $\Z^m/A$ and $\Z^n/B$ with the images of the standard basis as the generating sets. For $x \in \Z^m/A$ with an expression realizing the word norm $x = \sum_{k=1}^mz_ke_{k}$, we have 
\begin{align*}
|f(x)| = |\sum_{k=1}^mz_kf(e_{k})| \leq \sum_{k=1}^m|z_k||f(e_{k})| \leq (\max_k |f(e_k)|)\sum_{k=1}^m|z_k| = (\max_i |f(e_i)|)|x|.
\end{align*}
Hence this is Lipschitz. This completes the proof.
\end{proof}

However, the category $|{\rm Ab}|_{\Z}$ fails to be abelian. For example, let us consider two normed abelian groups $A$ and $B$ whose base groups are both $\Z$ with norms $|n|_A = |n|$ and $|n|_B= \begin{cases} 1 & n\neq 0 \\ 0 & n=0\end{cases}$. The identity on $\Z$ defines a Lipschitz homomorphism  $f : A \to B$, but we have
 \[
A \cong {\rm coker}({\rm ker} f) \not\cong {\rm ker}({\rm coker}f)\cong B,
\]
since $f$ is not an isomorphism, in fact, $f^{-1} : B \to A$ is not Lipschitz. To obtain such a standard isomorphism, we define the following condition.

\begin{df}
 For a morphism $f : A \to B$ in $|{\rm Ab}|_{\Z}$, we define
\[
|f^{-1}| = \inf \{c>0 \mid \forall a \in A, \inf_{a'\in  {\rm ker} f}|a+a'| \leq c|f(a)| \} \in [0, \infty],
\]
and we say $f$ is {\it strict} if $|f^{-1}|<\infty$. Equivalently, $f : A \to B$ is strict if the induced group isomorphism $A/{\rm ker}f \to {\rm im}f$ has a Lipschitz inverse.
\end{df}

\begin{lem}\label{freeabel}
Let $A \subset \Z^m$ and $B \subset \Z^n$ be subgroups. Then every group homomorphism $f : \Z^m/A \to \Z^n/B$ is a strict morphism in $|{\rm Ab}|_{\Z}$. 
\end{lem}
\begin{proof}
By Lemma \ref{zovera}, every group homomorphism $f : \Z^m/A \to \Z^n/B$ is Lipschitz. Now we show that the induced group isomorphism ${\rm coim} f \to {\rm im} f$ is an isomorphism in $|{\rm Ab}|_{\Z}$. 

Note that there exist subgroups $C \subset \Z^{m}$ and $D \subset \Z^{n}$ such that ${\rm ker} f = C/A$ and ${\rm im} f = D/B$, and hence ${\rm coim}f \cong \Z^m/C$ that preserves norms. Note that we have an isomorphism $D \cong \Z^{n'}$ in $|{\rm Ab}|_{\Z}$ for some $n'\geq0$. Indeed, after choosing a basis of the finite rank free abelian group $D$, the resulting group isomorphism $\Z^{n'}\to D$ and its inverse are Lipschitz by equivalence of norms on finite dimensional real vector spaces. Thus the induced group isomorphism ${\rm coim} f \to {\rm im} f$ is isomorphic to $\Z^{m}/C \to \Z^{n'}/B'$ for some  $B'$. Then Lemma \ref{zovera} implies that this is an isomorphism in $|{\rm Ab}|_{\Z}$. This completes the proof.
\end{proof}

\begin{lem}
For a morphism $f : A\to B$ in $|{\rm Ab}|_{\Z}$, the canonical homomorphism gives $ {\rm coker}({\rm ker} f)\cong {\rm ker}({\rm coker}f)$ in $|{\rm Ab}|_{\Z}$ if and only if $f$ is strict.
\end{lem}
\begin{proof}
Since we have $ {\rm coker}({\rm ker} f)\cong {\rm ker}({\rm coker}f)$ in ${\rm Ab}$, it holds in $|{\rm Ab}|_{\Z}$ if and only if there exist constants $C, c>0$ such that
\[
c|f(a)| \leq \inf_{a' \in {\rm ker}f}|a+a'| \leq C|f(a)|,
\]
for all $a \in A$. The first inequality holds since $f$ is  Lipschitz. Also, the second inequality holds if and only if  $f$ is strict. This completes the proof.
\end{proof}

\begin{rem}
Even if we restrict morphism class to strict morphisms, the resulting category is not abelian. In fact, it is not even a category since strictness is not preserved by compositions. For example, the composition $p\circ i : A \to A\oplus B \to B$, where $A$ and $B$ are as above and $i(n) = (n, n), p(m, n) = n$, is not strict although $i$ and $p$ are so.
\end{rem}

\subsubsection{The restricted products of normed abelian groups}

Next we define the restricted product of normed abelian groups.

\begin{df}\label{restrictedproduct}
\begin{enumerate}
\item Let $\{A_i\}_{i \in I}$ be a family of normed abelian groups. We define its {\it restricted product} by
\[
\overset{b}{\prod}_iA_i = \{(a_i)_i \in \prod_i A_i \mid \sup_i|a_i| < \infty\},
\]
with the norm defined by $|(a_i)_i | = \sup_i|a_i| $.
\item Let $\{f_i : A_i \to B_i\}_{i \in I}$ be a family of morphisms  in $|{\rm Ab}|_{\Z}$. It is {\it uniformly Lipschitz} if $\sup_i |f_i| < \infty$. If it additionally satisfies that $\sup_i |f^{-1}_i| < \infty$, we say it is {\it uniformly strict}.
\end{enumerate}
\end{df}

\begin{lem}\label{defrestprod}
Let $\{A_i\}_{i \in I}$ be a family of normed abelian groups. Then the restricted product $\overset{b}{\prod}_iA_i$ satisfies the following.
\begin{enumerate}
\item The family of projections $\{\pi_i\}_{i\in I}$ defined by
\[
\pi_i : \overset{b}{\prod}_iA_i \to A_i ; (a_i)_i \mapsto a_i
\]
is uniformly Lipschitz.
\item For every uniformly Lipschitz family of morphisms $\{f_i : X \to A_i\}_i$ in $|{\rm Ab}|_{\Z}$, there exists a unique morphism  $\varphi : X \to \overset{b}{\prod}_iA_i$ with $f_i = \pi_i\circ \varphi$ for all $i\in I$.
\end{enumerate}
\end{lem}
\begin{proof}
\begin{enumerate}
\item It is clear from $|\pi_i| \leq 1$ for all $i$.
\item We can define a group homomorphism by $\varphi(x) = (f_i(x))_i$ which is Lipschitz since $f_i$'s are uniformly Lipschitz. It is clear that $f_i = \pi_i\circ \varphi$ for all $i\in I$. The uniqueness is also clear.
\end{enumerate}
\end{proof}

\subsubsection{(AB$4^\ast$)-type exactness and a useful criterion}

\begin{lem}\label{kercoker}
Let $\{f_i : A_i \to B_i\}_{i \in I}$ be a uniformly Lipschitz family in $|{\rm Ab}|_{\Z}$. Let 
\[
\prod_if_i : \overset{b}{\prod}_iA_i \to \overset{b}{\prod}_iB_i
\]
be the morphism induced by Lemma \ref{defrestprod}. Then we have ${\rm ker}\prod_if_i = \overset{b}{\prod}_i{\rm ker}f_i$. Also, if $\{f_i \}_{i \in I}$ is uniformly strict, then we have ${\rm coker}\prod_if_i \cong \overset{b}{\prod}_i{\rm coker}f_i$, which preserves the norm.
\end{lem}
\begin{proof}
We have 
\[
{\rm ker}\prod_if_i = \{(a_i)_i \in \overset{b}{\prod}_iA_i  \mid \forall i, f_i(a_i) = 0 \} =\overset{b}{\prod}_i{\rm ker}f_i.
\]
Since both groups ${\rm ker}\prod_if_i$ and $\overset{b}{\prod}_i{\rm ker}f_i$ are subgroups of $\overset{b}{\prod}_iA_i $, the norms coincide. Also, the uniformly Lipschitz family $\{B_i \to {\rm coker}f_i\}_i$ induces a morphism $\overset{b}{\prod}_iB_i \to \overset{b}{\prod}_i{\rm coker}f_i$, which induces a morphism $\varphi : {\rm coker}\overset{b}{\prod}_if_i \to \overset{b}{\prod}_i{\rm coker}f_i$. The morphism $\varphi$ is obviously surjective: let $([b_i])_i \in \overset{b}{\prod}_i{\rm coker}f_i$. Then we can choose $b_i' \in B_i$ such that $[b_i] = [b_i'] \in {\rm coker}f_i$ and $(b_i')_i \in \overset{b}{\prod}_iB_i$, which implies $\varphi([(b_i')_i]) = ([b_i])_i $. Also, we can show the injectivity of $\varphi$. Note that, for $(a_i)_i \in \prod_iA_i$ with $(f_i(a_i))_i \in \overset{b}{\prod}_iB_i$, we can choose $(a_i')_i \in \overset{b}{\prod}_iA_i$ such that $f_i(a_i) = f_i(a_i')$ for all $i$ by the uniform strictness. In fact, for each $i$, we can choose $a_i' = a_i + k_i$ with $f(k_i) = 0$ such that $|a_i'|\leq c|f(a_i)|$, where $c$ is the constant associated with the uniform strictness of $\{f_i\}_i$. Hence we obtain
\begin{align*}
\exists (a_i)_i \in\overset{b}{\prod}_iA_i\text{ s.t. } (f(a_i))_i = (b_i-b'_i)_i  \Leftrightarrow \exists (a_i)_i \in \prod_iA_i \text{ s.t. }(f(a_i))_i = (b_i-b'_i)_i,
\end{align*}
which implies the injectivity of $\varphi$. Finally, note that we have $|\varphi([(b_i)_i])| \leq |[(b_i)_i]|$. Suppose that this is a strict inequality. Then, since the norms take values in $\Z_{\geq0}$, we can take $(b_i)_i \in \prod_iB_i$ and $(b_i')_i \in \overset{b}{\prod}_iB_i$ such that $\sup_i|b_i| < \sup_i|b'_i| = |[(b'_i)_i]|$ with $b_i-b_i' = f_i(a_i)$ for some $a_i \in A_i$ and for all $i$. We can suppose $(a_i)_i \in \overset{b}{\prod}_iA_i$ by the uniform strictness, and hence we obtain $[(b_i)_i] = [(b'_i)_i] \in {\rm coker}\overset{b}{\prod}_if_i$. However, we have $\sup_i|b'_i| = |[(b'_i)_i]|\leq \sup_i|b_i|$, which is a contradiction. This completes the proof.
\end{proof}
Note that, by Lemma \ref{kercoker}, we have
\begin{align*}
{\rm coim}\prod_if_i &=  {\rm coker}({\rm ker} \prod_if_i)  = {\rm coker}({\prod}_i{\rm ker} f_i) = \overset{b}{\prod}_i{\rm coker}({\rm ker} f_i) = \overset{b}{\prod}_i {\rm coim}f_i,\\
{\rm im}\prod_if_i &=  {\rm ker}({\rm coker} \prod_if_i)  = {\rm ker}(\prod_i{\rm coker} f_i) = \overset{b}{\prod}_i{\rm ker}({\rm coker} f_i) = \overset{b}{\prod}_i {\rm im}f_i.
\end{align*}
The following is an (AB$4^\ast$) -type exactness property for restricted products in $|{\rm Ab}|_{\Z}$.
\begin{lem}\label{normedab4}
Let us consider a family of short exact sequences
\[
0 \to A_i \xrightarrow{f_i} B_i \xrightarrow{g_i}C_i \to 0,
\]
in $|{\rm Ab}|_{\Z}$, where $\{f_i\}_i$ and $\{g_i\}_i$ are uniformly strict morphisms. Then the restricted product
\[
0 \to \overset{b}{\prod}_iA_i \xrightarrow{ \prod_if_i} \overset{b}{\prod}_iB_i \xrightarrow{ \prod_ig_i} \overset{b}{\prod}_iC_i\to 0,
\]
is exact.
\end{lem}
\begin{proof}
By Lemma \ref{kercoker}, we have ${\rm ker}\prod_if_i = \overset{b}{\prod}_i {\rm ker}f_i = 0$, and ${\rm coker}\prod_ig_i = \overset{b}{\prod}_i {\rm coker}g_i = 0$. Also, we have
\[
{\rm ker}\prod_ig_i = \overset{b}{\prod}_i {\rm ker}g_i = \overset{b}{\prod}_i {\rm im}f_i = {\rm im}\prod_if_i. 
\]
This completes the proof.
\end{proof}

\begin{lem}\label{prodchain}
Let $\{(C^i_\ast, \del^i_\ast)\}_{i\in I}$ be a family of chain complexes in  $|{\rm Ab}|_{\Z}$. Suppose that $\{\del^i_n\}_i$ is uniformly strict for all $n$. Then we have a chain complex 
\[
(\overset{b}{\prod}_iC^i_\ast, \prod_i\del^i_\ast),
\]
whose homology is isomorphic to $\overset{b}{\prod}_iH_\ast(C^i_\ast)$. Here we define homology group $H_n(C_\ast)$ as a normed abelian group ${\rm ker}\del_n/{\rm im}\del_{n+1}$.
\end{lem}
\begin{proof}
This is immediate from Lemma \ref{normedab4}.
\end{proof}

The following is obvious but gives a useful criterion for a family of morphisms to be uniformly strict.

\begin{lem}\label{fag}
Let $\{f_i:A_i\to B_i\}_{i\in I}$ be a family of strict morphisms in $|{\rm Ab}|_{\Z}$. Suppose that, up to norm-preserving isomorphisms of the source and target, only finitely many morphisms occur. Then $\{f_i\}_i$ is uniformly strict.
 \end{lem}

\subsection{Conventions for graphs}

\subsubsection{Graphs as metric spaces}
\begin{df}\label{defgraph}
\begin{enumerate}
\item {\it A graph} $X$ consists of a set $V(X)$, called the vertex set, together with a reflexive and symmetric relation $E(X)\subset V(X)\times V(X)$, called the edge relation. We identify $X$ with $V(X)$.

\item For graphs $X$ and $Y$, {\it a graph homomorphism} from $X \to Y$ is a map $f : V(X)\to V(Y)$ that induces a map of relations $E(X) \to E(Y)$. We write it $f : X \to Y$ by abbreviation. 
\item A graph is {\it connected} if the transitive closure of the relation $E(X)$ has just one equivalence class.
\item We regard a connected graph $X$ as a metric space equipped with {\it the path metric}, namely we define a distance function on $V(X)$ by
\[
d_X(x, x') = \min \{n \mid \exists x_0, \dots, x_n \text{ s.t. } x=x_0, x'=x_n, x_i\sim x_{i+1} \text{ for } 0\leq i <n\}.
\]
\end{enumerate}
When there is no ambiguity, we denote the distance function by $d$.
\end{df}

For an arbitrary graph $X$, we use the extended path metric, allowing the value $\infty$ when two vertices lie in different connected components. In the following, unless otherwise stated, the ambient graphs are assumed to be connected.

\begin{lem}
For graphs $X$ and $Y$, a map $f : X \to Y$ is a graph homomorphism if and only if it is a 1-Lipschitz map with respect to the path metric, namely we have $d_Y(f(x), f(x')) \leq d_X(x, x')$ for all $x, x' \in X$.
\end{lem}

In what follows, we mainly use the metric-space terminology and notation.

\subsubsection{The projecting decomposition}\label{projectingdecomp}

\begin{df}\label{bddgeom}
\begin{enumerate}
\item For a metric space $(X, d)$ and $r \geq0$, we define {\it the closed ball of radius $r$ centered at} $x$  by
\[
B_X(x, r) = \{x' \in X \mid d(x, x')\leq r\},
\]
for all $x \in X$. When there is no ambiguity, we denote it by $B(x, r)$.
\item A graph $X$ is {\it locally finite} if $B(x, r)$ is a finite set for all $x \in X$ and all $r\geq0$.
\item A graph $X$ has {\it bounded geometry} if it is locally finite and 
\[
\sup_{x \in X}\#B(x, r) < \infty
\]
for all $r\geq0$.
\end{enumerate}
\end{df}

\begin{df}
Let $X$ be a graph.
\begin{enumerate}
\item A graph $U$ is a {\it subgraph of} $X$ if $V(U)\subset V(X)$ and $E(U) \subset E(X)$. We write $U \subset X$ in this case.
\item A subgraph $U \subset X$ is {\it convex} if $d_U = d_X|_{U}$.
\item For $U \subset X$, we say $X$ {\it projects to} $U$ if there is a map $\pi : X \to U$ satisfying the following: for every $x \in X$ and every $u \in U$, we have 
\[
d_X(x, u) = d_X(x, \pi(x)) + d_X(\pi(x), u).
\]
We call such a $\pi$ {\it a projecting map}.
\item For $G, H\subset X$, subgraphs $G\cap H$ and $G\cup H$ are defined by taking $\cap$ and $\cup$, respectively,  of their vertex sets and edge sets. 
\item For $G, H \subset X$, we say $(X; G, H)$ is {\it a projecting decomposition} if  the following conditions hold :
\begin{itemize}
\item $G\cup H = X$,
\item $G\cap H$ is convex in $X$,
\item $G$ projects to $G\cap H$.
\end{itemize}
\end{enumerate}
\end{df}

\begin{lem}\label{pilip}
Let $U \subset X$ be a subgraph. Suppose that $X$ projects to $U$, and let $\pi : X \to U$ be a projecting map. Then we have the following.
\begin{enumerate}
\item $\pi|_U = {\rm id}_U$,
\item $\pi$ is $1$-Lipschitz.
\end{enumerate}
\end{lem}
\begin{proof}
\begin{enumerate}
\item For $u \in U$, we have $0 = d(u, u) = d(u, \pi(u)) + d(\pi(u), u)$, which implies $\pi(u) = u$.
\item For $x, x' \in X$, we have
\begin{align*}
d(\pi(x), x) + d(x, x') &\geq d(\pi(x), x') = d(\pi(x'), x') + d(\pi(x), \pi(x')) ,\\
d(\pi(x'), x') + d(x', x) &\geq d(\pi(x'), x) = d(\pi(x), x) + d(\pi(x'), \pi(x)) ,
\end{align*}
by the property of $\pi$. Adding these two inequalities, we obtain $d(x, x') \geq d(\pi(x), \pi(x'))$, which implies that $\pi$ is 1-Lipschitz.

\end{enumerate}
\end{proof}

\begin{rem}
The symmetry of the distance function is essential in Lemma \ref{pilip} (2); the statement need not hold for asymmetric distance functions.
\end{rem}
We recall the following lemma of Leinster \cite[Lemma~4.3]{L}, which will be useful later.

\begin{lem}[{\cite{L}, Lemma 4.3}] \label{Gconv}
Let $G, H \subset X$ with $G \cup H = X$. If $G\cap H$ is convex in $X$, then so are $G$ and $H$ in $X$.
\end{lem}

\section{Functions on tuples}\label{functionsontuples}
In this section, we encode possibly infinite formal sums of tuples as integer-valued functions on tuples.

\subsection{Conventions}
\begin{df}
 Let $X$ be a graph, and $n\geq 0$.
\begin{enumerate}
\item For a tuple $(x_0, \dots, x_n) \in X^{n+1}$, we define its {\it diameter} by ${\rm diam} (x_0, \dots, x_n) = \max_{i,j}d(x_i, x_j).$ Also, we define its {\it length} by $L(x_0, \dots, x_n) =\sum_{i=0}^{n-1}d(x_i, x_{i+1})$.
\item 
For a map $f : X^{n+1} \to \Z$, we define its support by
\[
{\rm Supp}f = f^{-1}(\Z\setminus\{0\}).
\]
We say that $f$ has {\it  controlled support} if we have
\[
\sup\{{\rm diam} (x_0, \dots, x_n) \mid (x_0, \dots, x_n) \in {\rm Supp}f  \} < \infty.
\]
Here and below, we use the convention that $\sup \emptyset = 0$.
\item 
For a map $f : X^{n+1} \to \Z$, we define its norm by
\[
|f|  = \sup_{(x_0, \dots, x_n) \in X^{n+1}} |f(x_0, \dots, x_n)| \in [0, \infty].
\]
We say that $f$ is  {\it  bounded} if $|f|<\infty$.
\end{enumerate}
\end{df}

\subsection{Simplicial structures}
\begin{lem}\label{facemap}
Let $X$ be a locally finite graph, and $n\geq 1$. For a map $f : X^{n+1} \to \Z$ with controlled support and $0\leq i \leq n$, we define a map 
\[
d_if : X^n \to \Z ; (x_0, \dots, x_{n-1}) \mapsto \sum_{y \in X}f(x_0, \dots, x_{i-1}, y, x_i, \dots, x_{n-1}).
\]
Then $d_if$ is well-defined and has controlled support. Also, when $X$ has bounded geometry,  there exists a constant $C(n, D)>0$ with $|d_if| \leq C(n, D)|f|$, where 
\[
D = \sup_{(x_0, \dots, x_n)  \in {\rm Supp}f }{\rm diam}(x_0, \dots, x_n).
\]
In particular, $d_if $ is  bounded whenever $f$ is.
\end{lem}
\begin{proof}
Since $f$ has controlled support, all $y$ contributing nonzero terms to the sum
\[
 \sum_{y \in X}f(x_0, \dots, x_{i-1}, y, x_i, \dots, x_{n-1}),
 \]
 lie in a ball of uniformly bounded radius around any fixed entry of the tuple. Since $X$ is locally finite, there are only finitely many such points $y$ for each $(x_0, \dots, x_{n-1}) \in X^n$. Thus $d_if $ is well-defined.  When $X$ has bounded geometry, the number of such points $y$ is uniformly bounded. This implies $|d_if| \leq C(n, D)|f|$ for some constant $C(n, D)>0$. Also, since deleting an entry does not increase the diameter, $d_if$ has controlled support whenever $f$ does. This completes the proof.
\end{proof}

\begin{lem}
Let $X$ be a locally finite graph. For a map $f : X^{n+1} \to \Z$ and $0\leq i \leq n$, we define a map 
\[
s_if : X^{n+2} \to \Z ; (x_0, \dots, x_{n+1}) \mapsto \begin{cases}f(x_0, \dots, x_{i}, \hat{x}_{i+1}, \dots, x_{n+1}) & \text{if\ } x_i= x_{i+1}, \\ 0 & \text{otherwise}.\end{cases}
\]
Then $s_if $ has controlled support whenever $f$ does. Also, we have $|s_if| = |f|$, in particular, $s_if $ is bounded  whenever $f$ is.
\end{lem}
\begin{proof}
It is immediate.
\end{proof}

\begin{lem}
Let $X$ be a locally finite graph. For a map $f : X^{n+1} \to \Z$ with controlled support,  we have $d_id_j f = d_{j-1}d_i f$ for $i<j$, $s_is_j f = s_{j}s_{i-1}$ for $i>j$ and 
\[
d_is_j f = \begin{cases} s_{j-1}d_i f & i<j, \\ f &i=j, j+1, \\  s_jd_{i-1}f & j+1<i. \end{cases}
\]
\end{lem}
\begin{proof}
It is straightforward.
\end{proof}

\begin{df}
 Let $X$ be a locally finite graph. We define a simplicial abelian group $C_\bullet (X) $ by
  \[
 C_n(X) = \{f : X^{n+1} \to \Z \mid f \text{ has controlled support}\},
 \]
 for $n\geq0$ with face maps $d_i$ and degeneracy maps $s_i$. When $X$ has bounded geometry, we also define a simplicial subgroup $ C^{\rm uf}_\bullet(X)$ of $C_\bullet (X)$ by
 \[
 C^{\rm uf}_n(X) = \{f \in C_n(X) \mid f \text{ is  bounded}\}.
 \]
 \end{df}

\subsection{The length filtration}

\begin{df}\label{filteredchaincomplex}
 Let $X$ be a locally finite graph. For $\ell \geq 0$, we define a simplicial subgroup $ F_\ell C_\bullet(X) \subset C_\bullet(X) $ by
  \[
 F_\ell C_n(X) = \{f \in C_n(X) \mid \sup_{(x_0, \dots, x_n)  \in {\rm Supp}f }L(x_0, \dots, x_n) \leq \ell\}.
 \]
Deleting an entry does not increase the length of a tuple, while repeating an entry does not change it. Hence these subgroups are preserved by all face and degeneracy maps. We also define $F_\ell C^{\rm uf}_\bullet(X)$  analogously. We  also set $F_\ell C_\ast(X)=0$ and $F_\ell C_\ast^{\rm uf}(X)=0$ for $\ell<0$.
\end{df}

In the following, for a simplicial abelian group $A_\bullet$, we write $A_\ast$ when we consider the associated chain complex.

\begin{rem}\label{simplicialobject}
Since $L(x_0, \dots, x_n)\leq \ell$ implies ${\rm diam}(x_0, \dots, x_n)\leq \ell$, Lemma \ref{facemap} shows that the face maps on $F_\ell C_\bullet^{\rm uf}(X)$ are Lipschitz. Moreover, $|s_i f|=|f|$. Hence $F_\ell C_\bullet^{\rm uf}(X)$ is a simplicial object in $|{\rm Ab}|_{\Z}$.
\end{rem}

\begin{lem}\label{exhaustive}
We have 
\[
 \bigcup_{\ell \geq 0}F_\ell C_n(X)  = C_n(X),
\]
and
\[
 \bigcup_{\ell \geq 0}F_\ell C^{\rm uf}_n(X) = C^{\rm uf}_n(X).
\]
In particular, $F_\bullet$ defines filtrations on the simplicial abelian groups $ C_\bullet(X)$ and $ C^{\rm uf}_\bullet(X)$.
\end{lem}

\begin{proof}
It is immediate for $n=0$, and we suppose $n\geq1$. For $f \in C_n(X)$, we have 
\[
{\rm diam} (x_0, \dots, x_n)  < c,
\]
for all $(x_0, \dots, x_n) \in {\rm Supp}f$ and for some $c>0$. Then we have
\[
L(x_0, \dots, x_n) < nc,
\]
for all $(x_0, \dots, x_n) \in {\rm Supp}f$, which implies that $f \in F_{nc}C_n(X)$, and $f \in F_{nc}C^{\rm uf}_n(X)$ if $f \in C^{\rm uf}_n(X)$.  This completes the proof.
\end{proof}

\section{Magnitude, path and coarse homology}\label{MPC}
In this section, we define several chain complexes consisting of functions on tuples of locally finite graphs, together with their homology. We first recall {\it coarse} homology and {\it uniformly finite} homology, originally introduced by Roe \cite{Roe} and Block--Weinberger \cite{blockweinberger}, respectively. We adopt the definitions for locally finite graphs due to Nowak--Yu \cite{NY}. We then introduce Borel--Moore and uniformly finite versions of {\it magnitude} homology and {\it path} homology, extending the original theories of Hepworth--Willerton \cite{HW} and Grigor’yan--Lin--Muranov--Yau \cite{GLMY} from finite graphs to locally finite graphs. 

\subsection{Coarse homology}
The following definition agrees with the chain complexes and homology groups defined in Definitions 7.1.1 and 7.2.2 of \cite{NY}.
\begin{df}\label{chom}
 Let $X$ be a locally finite graph. We denote the homology of $ C_\ast(X)$ by $ H_\ast(X)$, and we call it {\it the coarse homology of} $X$. Also, when $X$ has bounded geometry, we denote the homology of $ C^{\rm uf}_\ast(X)$ by $ H^{\rm uf}_\ast(X)$, and we call it {\it the uniformly finite homology of} $X$. 
\end{df}

 
\begin{rem}\label{remdege}
Let $D_n(X)\subset C_n(X)$ be the subgroup consisting of functions supported on degenerate tuples, that is, tuples $(x_0, \dots, x_n)$ satisfying $x_{i-1}=x_i$ for some $i$. We define $C_\ast^\Delta(X)=C_\ast(X)/D_\ast(X)$. As a graded abelian group, $C_n^\Delta(X)$ is naturally identified with the subgroup of $C_n(X)$ consisting of functions that vanish on degenerate tuples. We equip $C_\ast^\Delta(X)$ with the quotient filtration. By the standard normalization theorem  for simplicial abelian groups, the projection $\pi : C_\ast(X) \to C^\Delta_\ast(X)$ is a chain homotopy equivalence. Moreover, this is {\it a filtered chain homotopy equivalence}, namely  $\pi$, its homotopy inverse, and the chain homotopy all preserve the filtration $F_\bullet$. In fact, in the normalization theorem, we construct a chain homotopy equivalence $C_\ast \simeq N_\ast C$ and an isomorphism $N_\ast C \cong C_\ast/D_\ast$. The chain homotopy equivalence is constructed from compositions and finite linear combinations of the face and degeneracy maps $d_i$ and $s_i$, all of which preserve the filtration. The above isomorphism also preserves the filtration. The same construction applies to $C^{\rm uf}_\ast(X)$, and we denote the corresponding chain complex by $(C^{\Delta {\rm uf}}_\ast(X), \del^{\Delta}_\ast) $ .
\end{rem}

\begin{rem}\label{BMhomology}
The chain complex $C_\ast(X)$ of Definition \ref{chom} is related to so-called {\it Borel--Moore homology} as follows. We first recall the relation between locally finite homology and Borel--Moore homology. For a countable locally finite simplicial complex $K$, let $C_\ast^{\rm lf}(K)$ be the locally finite simplicial chain complex, whose chains are locally finite formal sums of {\it oriented} simplices, where we fix a total order on the set of vertices. By \cite[Proposition~4.7]{HughesRanicki}, the homology $H_\ast\bigl(C_\ast^{\rm lf}(K)\bigr)$ agrees with the locally finite singular homology $H_\ast^{\rm lf}(|K|)$ of its geometric realization. This homology is also referred to as Borel--Moore homology; see, for example, \cite[Section~2.2]{Maxim}, where it is defined using the filtered colimit of locally finite simplicial chain complexes over triangulations.

We now relate this construction to $C_\ast(X)$. For a locally finite graph $X$, let $V_R(X)$ be the {\it Vietoris--Rips complex} at scale $R$, which is a countable locally finite simplicial complex. Let $C_\ast^{\rm ord, lf}(V_R(X))$ be a simplicial chain complex that consists of locally finite sums of  {\it ordered} simplices. Namely, its $n$-chains are locally finite formal sums of ordered tuples $(x_0, \ldots, x_n)$ of diameter at most $R$, allowing repetitions. The standard normalization and ordered-oriented comparison give a chain homotopy equivalence $C_\ast^{\rm ord, lf}(V_R(X))\simeq C_\ast^{\rm lf}(V_R(X))$. Indeed, the usual chain homotopies are carried simplexwise and therefore extend to locally finite sums. Now we have
\[
\varinjlim_{R\to\infty}
C_n^{\rm ord,lf}(V_R(X))
=
\left\{
f:X^{n+1}\to\mathbb Z
\ \middle|
\sup_{(x_0,\ldots,x_n)\in {\rm Supp}f}
{\rm diam}(x_0, \ldots, x_n)<\infty
\right\},
\]
which is exactly $C_n(X)$ of Definition \ref{chom}. Therefore we obtain $C_\ast(X) \cong\varinjlim_{R\to\infty}
C_\ast^{\rm ord, lf}(V_R(X))$, and, up to the standard ordered--oriented chain homotopy equivalence, $C_\ast(X)$ is the filtered colimit of the Borel--Moore simplicial chain complexes of the Vietoris--Rips complexes of $X$.
\end{rem}

\begin{rem}\label{Roedefinition}
As mentioned in \cite[Section~7.6]{NY}, the homology $H_\ast(X)$ defined in Definition \ref{chom} is a special case of the more general theory of coarse homology introduced by Roe; see \cite[Section~5.5]{Roe}. Given a locally finite homology theory $h_\ast^{\mathrm{lf}}$, Roe defines the corresponding coarse homology by $\varinjlim_i h_\ast^{\mathrm{lf}}(|N(\mathcal U_i)|)$, where $\{\mathcal U_i\}_i$ is an {\it anti-\v{C}ech system}; see also \cite[Section~7.6]{NY}. The resulting homology is independent of the choice of the anti-\v{C}ech system.

For a locally finite graph $X$, we may take an increasing sequence $R_i\to\infty$ and the indexed ball covers $\mathcal B_{R_i}=\{B(x,R_i)\}_{x\in X}$ as the anti-\v{C}ech system. Identifying the vertices of the nerve $N(\mathcal B_R)$ with the center of the balls, we have simplicial inclusions
 $ V_R(X)\subset N(\mathcal B_R)\subset V_{2R}(X)$. Hence, by cofinality, we have $\varinjlim_i h_\ast^{\mathrm{lf}}(|N(\mathcal B_{R_i})|)\cong\varinjlim_i h_\ast^{\mathrm{lf}}(|V_{R_i}(X)|)$. By taking $h_\ast^{\rm lf}$ to be Borel--Moore homology in the sense of Remark \ref{BMhomology}, and using the exactness of filtered colimits, we obtain
\[
\varinjlim_i H_\ast^{\mathrm{BM}}(|V_{R_i}(X)|)
\cong
H_\ast\left(
\varinjlim_i C_\ast^{\mathrm{lf}}(V_{R_i}(X))
\right).
\]
By Remark \ref{BMhomology}, up to the standard ordered--oriented chain homotopy equivalence, the chain complex on the right is $C_\ast(X)$. Thus $H_\ast(X)$ is precisely the coarse homology obtained by coarsening Borel--Moore homology.
\end{rem}

\subsection{Magnitude homology}
We next pass from the total complex to the associated graded complex of the length filtration.
\begin{df}\label{mag}
 Let $X$ be a locally finite graph. For $\ell \geq 0$, we denote the homology of a chain complex  $F_\ell C_\ast(X)/F_{\ell-1} C_\ast(X)$ by ${\rm MH}^{\rm BM}_{ \ast, \ell}(X)$, and we call it {\it the Borel--Moore magnitude homology of} $X$.  When $X$ has bounded geometry, we denote the homology of $F_\ell C^{\rm uf}_\ast(X)/F_{\ell-1} C^{\rm uf}_\ast(X)$ by ${\rm MH}^{\rm uf}_{\ast, \ell}(X)$, and we call it  {\it the uniformly finite magnitude homology of} $X$.  
 
 Explicitly,  they are the homology of the chain complex $({\rm MC}^{\rm BM}_{\ast, \ast}(X), \del^{\rm mag}_\ast)$ and  its subcomplex $({\rm MC}^{\rm uf}_{\ast, \ast}(X), \del^{\rm mag}_\ast)$ defined by 
  \[
{\rm MC}^{\rm BM}_{n, \ell}(X) = \{f \in C_n(X) \mid L(x_0, \dots, x_n) = \ell \text{ for } (x_0, \dots, x_n)  \in {\rm Supp}f \},
 \]
and
  \[
{\rm MC}^{\rm uf}_{n, \ell}(X) = \{f \in C^{\rm uf}_n(X) \mid L(x_0, \dots, x_n) = \ell \text{ for } (x_0, \dots, x_n)  \in {\rm Supp}f \},
 \]
 with
  \begin{align*}
 \del^{\rm mag}_{n, \ell} f(x_0, \dots, x_{n-1}) = \sum_{i=0}^n (-1)^i \sum_{x_{i-1}\leq y \leq x_i}f(x_0, \dots, x_{i-1}, y, x_i, \dots,  x_{n-1}).
 \end{align*}
 Here, for $i=0, n$, the conditions $x_{i-1}\leq y \leq x_i$ mean $y = x_0$ and $y=x_{n-1}$ respectively. Also, for $a, b, c \in X$, we write $a \leq b \leq c$ if $d(a, b) + d(b, c) = d(a, c)$. Note that the chain complex $({\rm MC}^{\rm uf}_{\ast, \ast}(X), \del^{\rm mag}_\ast)$ is only defined when $X$ has bounded geometry.
 \end{df}
 
 \begin{rem}\label{Delta}
As in Remark \ref{remdege}, we may equivalently use the definition of ${\rm MC}^{\rm BM}_{n, \ell}(X)$ and ${\rm MC}^{\rm uf}_{n, \ell}(X)$ by adding the condition
 \[
f(x_0, \dots, x_n) = 0 \text{ if } x_{i-1} = x_i \text{ for some }i,
\]
where the resulting chain complexes $({\rm MC}^{\Delta{\rm BM}}_{\ast, \ell}(X), \del^{\Delta{\rm mag}}_{\ast, \ell})$ and   $({\rm MC}^{\Delta{\rm uf}}_{\ast, \ell}(X), \del^{\Delta{\rm mag}}_{\ast, \ell})$ are chain homotopy equivalent to the original ones. It is immediate to verify that ${\rm MC}^{\Delta{\rm BM}}_{n, \ell}(X) = {\rm MC}^{\Delta{\rm uf}}_{n, \ell}(X) = 0$ for $n>\ell$, and hence we have ${\rm MH}^{\rm BM}_{n, n}(X) \cong {\rm ker}\del^{\Delta{\rm mag}}_{n, n}$ and similarly for ${\rm MH}^{\rm uf}_{n, n}(X)$. We also note that $F_\ell C_n^\Delta(X) = 0$ for $\ell < n$, $F_nC_n^\Delta(X) = {\rm MC}^{\Delta{\rm BM}}_{n, n}(X)$ and similarly for the uniformly finite version.
 \end{rem}
 
 \begin{rem}\label{ordinary}
 When $X$ is a finite graph, the chain complexes ${\rm MC}^{\rm BM}_\ast(X)$ and ${\rm MC}^{\rm uf}_\ast(X)$ coincide with the original definition of magnitude chain complex of $X$ due to Hepworth--Willerton \cite{HW}. Explicitly, the original magnitude chain complex $({\rm MC}_{\ast, \ell}(X), \del^{\rm mag}_{\ast, \ell})$ is defined as
 \[
{\rm MC}_{n, \ell}(X) = \Z\langle (x_0, \dots, x_n) \mid \sum_{i=0}^{n-1}d(x_i, x_{i+1}) = \ell \rangle,
\]
and
\[
 \del^{\rm mag}_{n, \ell} (x_0, \dots, x_{n}) = \sum_{\substack{0\leq i\leq n \\ x_{i-1}\leq x_i \leq x_{i+1}}} (-1)^i (x_0, \dots, \hat{x}_{i}, \dots,  x_{n}),
\]
Here we adopt the un-normalized version in the sense of Leinster--Shulman's Definition 5.7 in \cite{LS}, although Hepworth--Willerton defined it in a normalized form. The point is that the finite linear combination of tuples can be considered as a function on tuples with finite support. Thus Definition \ref{mag} extends magnitude homology to locally finite graphs possibly with infinitely many vertices. Hepworth--Willerton's definition also makes sense for an arbitrary graph $X$, using finite linear combinations of tuples. We denote its homology by ${\rm MH}_{\ast,\ast}(X)$ and refer to it as {\it the ordinary magnitude homology} in the following.
\end{rem}
For later use, we introduce the magnitude homology with fixed initial point here. For a finite graph $X$ and $x \in X$, we define
 \[
{\rm MC}^x_{n, \ell}(X) = \Z\langle (x_0, \dots, x_n) \mid x_0=x, \sum_{i=0}^{n-1}d(x_i, x_{i+1}) = \ell \rangle,
\]
and $\del^x_{n, \ell}$ as the restriction of $\del^{\rm mag}_{n, \ell}$, which together form a chain subcomplex ${\rm MC}^x_{\ast, \ell}(X)$ of ${\rm MC}_{\ast, \ell}(X)$ . In fact, it is a direct summand of $({\rm MC}_{\ast, \ell}(X), \del^{\rm mag}_{\ast, \ell} )$ satisfying ${\rm MC}_{\ast, \ell}(X) = \oplus_x {\rm MC}^x_{n, \ell}(X) $. We denote its homology by ${\rm MH}^x_{n, \ell}(X)$. Although we defined it for finite graphs, we can also consider it as a definition for arbitrary graphs. In the case that  $X$ is locally finite,  every tuple of length $\ell$ with initial point $x$ is contained in the finite ball $B(x, \ell)$. Hence, for fixed $n$ and $\ell$, there are only finitely many such tuples, and  ${\rm MC}^x_{n, \ell}(X)$ has finite rank. 

\begin{prp}\label{prod}
For a locally finite graph $X$, we have 
\[
{\rm MH}^{\rm BM}_{\ast, \ell}(X) \cong \prod_{x \in X} {\rm MH}^{x}_{\ast, \ell}(X).
\]
Also, when $X$ has bounded geometry, we have
\[
{\rm MH}^{\rm uf}_{\ast, \ell}(X) \cong \overset{b}{\prod}_{x \in X} {\rm MH}^{x}_{\ast, \ell}(X).
\]
Here we consider ${\rm MC}^{\rm uf}_{\ast, \ell}$  as  a chain complex in $|{\rm Ab}_{\Z}|$ as in Remark \ref{simplicialobject} and equip the homology on the left-hand side with the quotient norm explained in Section \ref{prelimabelian}. Also, each ${\rm MC}_{n,\ell}^x(X)$ is equipped with the $\ell^\infty$-norm with respect to the basis of tuples, and its homology is equipped with the quotient norm as in Section \ref{prelimabelian}. Moreover, the following diagram is commutative:
\begin{equation*}
\xymatrix{
{\rm MH}^{\rm uf}_{\ast, \ell}(X) \ar[d]^{\iota_\ast} \ar[r]^{\cong} & \overset{b}{\prod}_{x \in X} {\rm MH}^{x}_{\ast, \ell}(X) \ar@{^{(}-_>}[d] \\
{\rm MH}^{\rm BM}_{\ast, \ell}(X) \ar[r]^{\cong} & \prod_{x \in X} {\rm MH}^{x}_{\ast, \ell}(X),
}
\end{equation*}
where $\iota_\ast$ is induced from the inclusion ${\rm MC}^{\rm uf}_{\ast, \ell}(X) \to {\rm MC}^{\rm BM}_{\ast, \ell}(X)$.
\end{prp}
\begin{proof}
Note that  we have an isomorphism
\[
{\rm MC}^{{\rm BM}}_{\ast, \ell}(X) \cong  \prod_{x \in X} {\rm MC}^{x}_{\ast, \ell}(X).
\]
Hence Lemma \ref{ab4} implies the isomorphism for Borel--Moore magnitude homology. Also, when $X$ has bounded geometry, the above isomorphism induces 
\[
{\rm MC}^{{\rm uf}}_{\ast, \ell}(X) \cong  \overset{b}{\prod}_{x \in X} {\rm MC}^{x}_{\ast, \ell}(X).
\]
For fixed $n$ and $\ell$, bounded geometry implies that the ranks of ${\rm MC}_{n,\ell}^x(X)$ and ${\rm MC}_{n-1,\ell}^x(X)$ are uniformly bounded. Indeed, every tuple contributing to ${\rm MC}^x_{n, \ell}(X)$ is contained in $B(x,\ell)$, whose cardinality is uniformly bounded in $x$ by bounded geometry. Moreover, with respect to the bases given by tuples, the matrices representing $\partial_{n,\ell}^x$ have entries in $\{0,\pm1\}$. Since the ranks are uniformly bounded and the matrix entries belong to a fixed finite set, up to permutations of the tuple bases only finitely many morphisms $\del^x_{n,\ell}$ occur. Note that such permutations preserve the $\ell^\infty$-norms, and each $\partial^x_{n,\ell}$ is strict by Lemma \ref{freeabel}. Therefore Lemma \ref{fag} implies that $\{\del^x_{n, \ell}\}_{x\in X}$ is uniformly strict, and  then Lemma \ref{prodchain} implies the desired isomorphism. The diagram above is commutative since it is already commutative at the chain level. This completes the proof.
\end{proof}

We recall here the definition of {\it diagonality} of a graph due to Hepworth--Willerton \cite[Definition~7.1]{HW}. A graph is {\it diagonal} if we have 
\[
{\rm MH}_{n, \ell}(X) = 0 \text{ for } \ell \neq n.
\]
Now we define Borel--Moore and uniformly finite versions of diagonality.
\begin{df}
Let $X$ be a locally finite graph. We say ${\rm MH}^{\rm BM}_{\ast, \ast}(X)$ is {\it diagonal} if it satisfies
\[
{\rm MH}^{\rm BM}_{n, \ell}(X) = 0 \text{ for } \ell \neq n.
\]
Also, when $X$ has bounded geometry, we say ${\rm MH}^{\rm uf}_{\ast, \ast}(X)$ is {\it diagonal} if it satisfies
\[
{\rm MH}^{\rm uf}_{n, \ell}(X) = 0 \text{ for } \ell \neq n.
\]
\end{df}

\begin{prp}\label{diagonality}
Let $X$ be a locally finite graph. Then ${\rm MH}^{\rm BM}_{\ast, \ast}(X)$  is  diagonal if and only if $X$ is diagonal. Also, when $X$ has bounded geometry, ${\rm MH}^{\rm uf}_{\ast, \ast}(X)$ is  diagonal if and only if $X$ is diagonal.
\end{prp}
\begin{proof}
Since we have ${\rm MH}_{\ast, \ast}(X) \cong \oplus_x {\rm MH}^x_{\ast, \ast}(X)$,  ${\rm MH}^{\rm BM}_{\ast, \ell}(X) \cong \prod_{x \in X} {\rm MH}^{x}_{\ast, \ell}(X)$, and $ {\rm MH}^{\rm uf}_{\ast, \ell}(X) \cong \overset{b}{\prod}_{x \in X} {\rm MH}^{x}_{\ast, \ell}(X)$ when $X$ has bounded geometry, by Proposition \ref{prod}, all of them are diagonal if and only if so is each component ${\rm MH}^{x}_{\ast, \ast}(X)$. This completes the proof. 
\end{proof}

\subsection{Path homology}
We now extract the diagonal part of the normalized magnitude complex and relate it to path homology. In what follows, we regard $\mathrm{MC}^{\Delta\mathrm{BM}}_{n,\ell}(X)$ and $\mathrm{MC}^{\Delta\mathrm{uf}}_{n,\ell}(X)$ as subgroups of $C^\Delta_n(X)$ and $C^{\Delta\mathrm{uf}}_n(X)$, respectively, via the identification of chains with functions on tuples as in Remarks \ref{remdege} and \ref{Delta}. 

\begin{df}\label{path}
Let $X$ be a locally finite graph. 
\begin{enumerate}
\item For $n\geq 0$, we define 
\begin{align*}
A_n^{\rm BM}(X) &= {\rm MC}^{\Delta{\rm BM}}_{n, n}(X) \\
&=\{f \in C^{\Delta}_n(X) \mid f(x_0, \dots, x_n) = 0 \text{ if } d(x_i, x_{i+1}) \neq 1 \text{ for some }i\}.
\end{align*}
Also, when $X$ has bounded geometry, we define
\begin{align*}
A_n^{\rm uf}(X) &= {\rm MC}^{\Delta{\rm uf}}_{n, n}(X) \\
&=\{f \in C_n^{\Delta{\rm uf}}(X) \mid f(x_0, \dots, x_n) = 0 \text{ if } d(x_i, x_{i+1}) \neq 1 \text{ for some }i\}.
\end{align*}
\item We define a chain complex $(\Omega^{\rm BM}_\ast(X), \del^{\rm path}_{\ast})$  by $\Omega^{\rm BM}_0(X) = A_{0}^{\rm BM}(X)$,
\[
\Omega^{\rm BM}_n(X) = A_{n}^{\rm BM}(X)\cap(\del^{\Delta}_n)^{-1}A_{n-1}^{\rm BM}(X), 
\]
for $n\geq1$, and
\[
\del^{\rm path}_n = (\del^\Delta_n)|_{\Omega^{\rm BM}_n(X) },
\]
where $\del^{\Delta}_n : C^\Delta_n(X) \to C^\Delta_{n-1}(X)$ is defined in Remark \ref{remdege}. Here we use the identification ${\rm MC}_{n, n}^{\Delta{\rm BM}} \subset C^\Delta_n(X)$. Also, when $X$ has bounded geometry, we similarly define a chain complex $(\Omega^{\rm uf}_\ast(X), \del^{\rm path}_{\ast})$ by $\Omega^{\rm uf}_0(X) = A_{0}^{\rm uf}(X)$,
\[
 \Omega^{\rm uf}_n(X) = A_{n}^{\rm uf}(X)\cap(\del^{\Delta}_n )^{-1}A_{n-1}^{\rm uf}(X), 
 \]
for $n\geq1$, and
 \[
\del^{\rm path}_n = (\del^\Delta_n)|_{\Omega^{\rm uf}_n(X)}.
\]
Since $(\partial^\Delta)^2=0$, we have $\partial_n^\Delta(\Omega_n^{\rm BM})\subset \Omega_{n-1}^{\rm BM}$, and similarly for the uniformly finite version. We denote their homology groups by ${\rm PH}^{\rm BM}_n(X)$ and ${\rm PH}^{\rm uf}_n(X)$, respectively, and we call them {\it the Borel--Moore  path homology } and {\it the uniformly finite path homology} of $X$.
\end{enumerate}
\end{df}

 \begin{rem}
 By the same discussion in Remark \ref{ordinary}, when $X$ is a finite graph considered as a digraph with bi-directed edges, the chain complexes $(\Omega^{\rm BM}_\ast(X), \del^{\rm path}_{\ast})$ and $(\Omega^{\rm uf}_\ast(X), \del^{\rm path}_{\ast})$ coincide with the original definition of path chain complex of $X$ due to Grigor'yan--Lin--Muranov--Yau \cite{GLMY}. Hence Definition \ref{path}  extends path homology for undirected finite graphs to locally finite graphs possibly with infinitely many vertices.
  \end{rem}
  
\begin{lem}\label{omemh}
For a locally finite graph $X$, we have $\Omega^{\rm BM}_n(X) \cong {\rm MH}^{\rm BM}_{n, n}(X)$ and $\Omega^{\rm uf}_n(X)  \cong {\rm MH}^{\rm uf}_{n, n}(X)$, where the latter is defined when $X$ has bounded geometry.
\end{lem}
\begin{proof}
It is immediate for $n=0$, and suppose $n\geq1$. By Remark \ref{Delta}, we have ${\rm MH}^{\rm BM}_{n, n}(X) \cong {\rm ker}\del^{\Delta{\rm mag}}_{n, n}$ and $F_nC^\Delta_n(X) = {\rm MC}^{\Delta{\rm BM}}_{n, n}(X)$. Hence, from the definitions, we obtain
\begin{align*}
{\rm MH}^{\rm BM}_{n, n}(X) &\cong {\rm ker}\del^{\Delta{\rm mag}}_{n, n} \\
&=   {\rm MC}^{\Delta{\rm BM}}_{n, n}(X)\cap(\del^{\Delta}_n)^{-1}(F_{n-1}C^\Delta_{n-1}(X)) \\
&= A^{\rm BM}_n(X) \cap (\del^{\Delta}_n)^{-1}A_{n-1}^{\rm BM}(X) \\
&= \Omega^{\rm BM}_n(X).
\end{align*}
It is similar for $\Omega^{\rm uf}_n(X) $. This completes the proof.
\end{proof}

\section{Borel--Moore and uniformly finite magnitude--path spectral sequences}\label{BMUFMPSS}
In this section, we construct the {\it Borel--Moore} and {\it uniformly finite magnitude--path spectral sequences} associated with the length filtrations on $C_\ast(X)$ and $C_\ast^{\rm uf}(X)$. We first recall in Section \ref{spectralseq} the spectral-sequence conventions used in this paper, with particular emphasis on convergence for filtrations that are bounded below and exhaustive but not necessarily bounded. In Section \ref{sec:5.2}, we identify the $E^1$-pages with Borel--Moore and uniformly finite magnitude homology, the $p$-axis of the $E^2$-pages with the corresponding path homology, and the $E^\infty$-pages with the associated graded groups of coarse and uniformly finite homology. We then study functoriality and invariance properties of these spectral sequences in Section \ref{sec:5.3}. Finally, Section \ref{sec:5.4} gives an example in which the spectral sequence  converges although a term does not stabilize at any finite page.

\subsection{Spectral sequences associated to filtered chain complexes}\label{spectralseq}
We recall the spectral-sequence conventions used in this paper, with particular emphasis on convergence. The length filtrations considered below are exhaustive and bounded below, but not bounded, so the terms of the spectral sequence need not stabilize at any finite page. See \cite{Wei} for details. 

A {\it spectral sequence} in an Abelian category is a family of objects $\{E^r_{p, q}\}$ for $p, q \in \Z$ and $r\geq 0$ with morphisms $d^r_{p, q} : E^r_{p, q} \to E^r_{p-r, q+r-1}$ such that
\begin{itemize}
\item $d^r_{\ast, \ast}$'s are differentials, namely $d^r_{p-r, q+r-1}\circ d^r_{p, q} = 0$,
\item $E^{r+1}_{\ast, \ast}$'s are isomorphic to the homology of $E^r_{\ast, \ast}$, namely $E^{r+1}_{p, q} \cong {\rm ker}d^r_{p, q}/{\rm im}d^r_{p+r, q-r+1}$.
\end{itemize}

Given a filtered chain complex $F_\bullet C_\ast$, where the filtration is indexed over $\Z$, we can construct a spectral sequence as follows (\cite{Wei} Section 5.5). For $r \in \Z$, we define
\begin{align*}
Z^r_{p,q} &= F_pC_{p+q}\cap \del_{p+q}^{-1}(F_{p-r}C_{p+q-1}), \\
B^r_{p, q} &= F_pC_{p+q}\cap \del_{p+q+1}F_{p+r}C_{p+q+1},
\end{align*}
and
\[
E^r_{p, q} = \frac{Z^r_{p, q}}{Z^{r-1}_{p-1, q+1} + B^{r-1}_{p, q}}.
\]
Note that we have 
\[
E^r_{p, q} =F_pC_{p+q}/F_{p-1}C_{p+q} \text{ for } r\leq 0,
\]
and 
\[
E^1_{p, q} = H_{p+q}(F_pC_{\ast}/F_{p-1}C_{\ast}).
\]
Also, the differential on $C_\ast$ induces $d^r_{p, q} : E^r_{p, q} \to E^r_{p-r, q+r-1}$ for all $p, q, r$, which satisfies the above two conditions. We note that a filtered chain map $f : F_\bullet C_\ast \to F_\bullet D_\ast$ induces homomorphisms $f^r_{p, q} : E^r_{p, q}(C_\ast) \to E^r_{p, q}(D_\ast)$ satisfying $f^r_{p-r, q+r-1}\circ d^r_{p, q} = d^r_{p, q}\circ f^r_{p, q}$. Therefore, if there exists $r$ such that $f^r_{p, q}$'s are isomorphisms for all $p, q$, then so are $f^R_{p, q}$'s for $R\geq r$. In particular, if $f$ is a filtered chain homotopy equivalence, then the induced homomorphisms $f^r_{p, q}$'s are isomorphisms for all $p, q$ and all $r\geq1$. 

A virtue of a spectral sequence becomes apparent when we consider its {\it convergence}. A typical situation is that the filtration $F_\bullet$ is {\it bounded}, namely, for each $n \in \Z$, there exist integers $s, t$ such that $F_sC_n = 0$ and $F_tC_n = C_n$. In this case, the directions of the differentials $d^r$ force each $E^r_{p, q}$ to stabilize for sufficiently large $r$, and we denote the stabilized term by $E^\infty_{p, q}$. Moreover, the homology $H_n(C_\ast)$ of the original chain complex $C_\ast$ can be determined up to {\it extension problems} from the terms $E^\infty_{p, n-p}$.

In the case of the filtered chain complexes defined in Definition \ref{filteredchaincomplex}, however, the filtration is not bounded. Instead, it is {\it bounded below}, namely we have $F_pC_n=0$ for $p<0$. In this case, even though the terms $E^r_{p,q}$ need not stabilize at any finite stage, we can still discuss the convergence of the spectral sequence as follows. For a filtered chain complex $F_\bullet C_\ast$, we define
\[
Z^{\infty}_{p, q} = \bigcap_{r=0}^\infty Z^r_{p, q}, B^{\infty}_{p, q} = \bigcup_{r=0}^{\infty} B^r_{p, q},
\]
and
\[
E^\infty_{p, q} = \frac{Z^\infty_{p, q}}{Z^\infty_{p-1, q+1} + B^\infty_{p, q}}.
\]
We also define a filtration $F_\bullet$ on $H_\ast(C_\ast)$ by
\[
F_pH_n(C_\ast) = \{\alpha \in H_n(C_\ast) \mid \exists z \in F_pC_n \text{ s.t. } \alpha = [z]\}.
\]
\begin{df}[\cite{Wei} 5.2.11]
 We say a spectral sequence $E^\ast_{\ast, \ast}$ associated with a bounded below filtered chain complex $C_\ast$ {\it  converges} if it satisfies $\cup_pF_pH_n(C_\ast) = H_n(C_\ast), \cap_pF_pH_n(C_\ast) = 0$ and 
\[
E^\infty_{p, q} \cong F_{p}H_{p+q}(C_\ast)/F_{p-1}H_{p+q}(C_\ast).
\]
\end{df}

A filtration $F_\bullet$ on a chain complex $C_\ast$ is called {\it exhaustive} if it satisfies $\cup_pF_pC_\ast = C_\ast$. The filtrations defined in Definition \ref{filteredchaincomplex} are exhaustive by Lemma \ref{exhaustive} and are bounded below, hence the corresponding spectral sequences converge by the following fact.

\begin{prp}[\cite{Wei} Theorem 5.5.1]\label{conv}
Let $E^\ast_{\ast, \ast}$  be a spectral sequence associated with a filtered chain complex $F_\bullet C_\ast$. If the filtration $F_\bullet$ is bounded below and exhaustive,  then $E^\ast_{\ast, \ast}$  converges.
\end{prp}

\begin{rem}\label{stability}
Even when a spectral sequence $E^\ast_{\ast, \ast}$  converges, each term  $E^r_{p, q}$ need not stabilize for sufficiently large $r$. For example, in the following filtered chain complex, we have $E^r_{0, 0} \not\cong 0$ for all $r$, while the spectral sequence converges and $E^\infty_{0, 0} \cong 0$ : 
\[
C_0 = \oplus_{k\geq 0} \Z x_k, C_1= \oplus_{k\geq 0} \Z y_k, C_n = 0 \text{ for } n\neq 0, 1,
\]
with $\del_1y_k = x_k$ and $\del_n = 0$ for $n \neq 0, 1$. We define a filtration $F_\bullet$ by
\[
F_pC_0 = C_0, F_pC_1 =  \oplus_{k= 0}^p \Z y_k,
\]
for all $p\geq 0$, and $F_pC_n = 0$ for $p<0$.
\end{rem}

Unlike the bounded case, the $E^\infty$-page does not in general determine $H_\ast(C_\ast)$ through a finite collection of extension problems. Nevertheless, it still allows us to compare the homologies of filtered chain complexes via the following comparison theorem.

\begin{prp}[\cite{Wei} Theorem 5.2.12]\label{comparison}
Let $F_\bullet C_\ast$ and $F_\bullet D_\ast$ be filtered chain complexes whose filtrations are bounded below and exhaustive. If a filtered chain map $f : F_\bullet C_\ast \to F_\bullet D_\ast$ induces isomorphisms of all $(p, q)$-terms $f^r_{p, q} : E^r_{p, q}(C_\ast) \cong E^r_{p, q}(D_\ast)$ for some $r\geq0$, then $f$ induces an isomorphism $f_\ast : H_\ast(C_\ast) \cong H_\ast(D_\ast)$.
\end{prp}

\subsection{Magnitude--path spectral sequences}\label{sec:5.2}
Now we apply the framework introduced in the previous subsection to construct the {\it Borel--Moore} and {\it uniformly finite magnitude--path spectral sequences}, which generalize the spectral sequence studied in \cite{A} from finite graphs to locally finite graphs.
\begin{df}
Let $X$ be a locally finite graph. We denote the spectral sequences associated to the filtered chain complexes $F_\bullet C_\ast(X)$ and $F_\bullet C^{\rm uf}_\ast(X)$ in Definition \ref{filteredchaincomplex}, where the latter is defined when $X$ has bounded geometry,  by $E^{{\rm BM}, r}_{p, q}(X)$ and $E^{{\rm uf}, r}_{p, q}(X)$, respectively.
\end{df}

\subsubsection{$E^1$-page as magnitude homology and  convergence to coarse homology}
The following proposition shows that the spectral sequences above  converge to the corresponding coarse homologies and that their $E^1$-terms are the Borel--Moore and uniformly finite magnitude homology groups, respectively. The latter statement is a generalization of  \cite[Proposition~7.8]{A}.

\begin{prp}\label{weaklyconverge}
For a locally finite graph $X$, the spectral sequences $E^{{\rm BM}, r}_{p, q}(X)$ and $E^{{\rm uf}, r}_{p, q}(X)$, where the latter is defined when $X$ has bounded geometry,  converge. Namely we have
\begin{align*}
E^{{\rm BM}, \infty}_{p, q}(X) &\cong F_pH_{p+q}(X)/F_{p-1}H_{p+q}(X),\\
E^{{\rm uf}, \infty}_{p, q}(X) &\cong F_pH^{\rm uf}_{p+q}(X)/F_{p-1}H^{\rm uf}_{p+q}(X),
\end{align*}
for $p\geq0$ and $p+q\ge0$. We also have
\begin{align*}
E^{{\rm BM}, 1}_{p, q}(X) &\cong {\rm MH}^{\rm BM}_{p+q, p}(X),\\
E^{{\rm uf}, 1}_{p, q}(X) &\cong {\rm MH}^{\rm uf}_{p+q, p}(X),
\end{align*}
for $p\geq0$ and $p+q\ge0$. 
\end{prp}
\begin{proof}
The convergence follows from Proposition \ref{conv}, since the length filtrations are bounded below by Definition \ref{filteredchaincomplex} and exhaustive by Lemma \ref{exhaustive}. The identifications of the $E^1$-terms follow immediately from Definition \ref{mag}. This completes the proof.
\end{proof}

\subsubsection{$E^2$-page is path homology}
We next show that the $E^2$ terms on the $p$-axis are exactly the path homology groups. This extends \cite[Proposition~6.11,~7.8]{A} from finite graphs to locally finite graphs. The only additional point is that our chains may have infinite support. However, by local finiteness and the controlled-support condition, for each fixed tuple only finitely many terms contribute to each face map, as shown in Lemma \ref{facemap}. Hence all boundary maps appearing in the finite-case argument remain well-defined, and the same chain-level argument applies to the Borel--Moore version. For the uniformly finite version, bounded geometry further guarantees that the number of contributing terms is uniformly bounded, so that bounded chains are preserved.

\begin{prp}\label{pathhomology}
For a  locally finite graph $X$, we have
\begin{align*}
E^{{\rm BM}, 2}_{p, 0}(X) &\cong {\rm PH}^{\rm BM}_{p}(X),\\
E^{{\rm uf}, 2}_{p, 0}(X) &\cong {\rm PH}^{\rm uf}_{p}(X),
\end{align*}
for all $p\geq0$, where the latter is defined when $X$ has bounded geometry.
\end{prp}
\begin{proof}
Since the canonical projection $C_\ast(X)\to C_\ast^\Delta(X)$ is a filtered chain homotopy equivalence as explained in Remark \ref{remdege}, it induces an isomorphism of the associated spectral sequences from the $E^1$-page. So we consider the spectral sequence  associated with $C^\Delta_{\ast}(X)$ in the following. By Lemma \ref{omemh} and Proposition \ref{weaklyconverge}, we have isomorphisms $E^{{\rm BM}, 1}_{p, 0}(X) \cong {\rm MH}^{\rm BM}_{p, p}(X) \cong \Omega^{{\rm BM}}_{p}(X)$. Under these isomorphisms, the differential $d^1_{p, 0}$ on the left term corresponds to $\partial^{\rm path}_p$ on the right term. In fact, for $p\geq1$, $f\in\Omega_p^{\mathrm{BM}}(X) \subset C_p^\Delta(X)$ corresponds to an element represented by the same chain
\[
f_{\rm mag} \in {\rm MH}^{\rm BM}_{p, p}(X)\cong {\rm ker}\del^{\Delta{\rm mag}}_p \subset {\rm MC}^{\Delta {\rm BM}}_{p, p}(X)=F_pC_p^\Delta(X),
\]
by Remark \ref{Delta}, and we have $\del^{\rm path}_p f = \del_p^\Delta f\in \Omega^{\rm BM}_{p-1}(X) \subset F_{p-1}C_{p-1}^\Delta(X)$. On the other hand, $d^1(f_{\rm mag})$ is represented by $\del_p^\Delta f \in F_{p-1}C_{p-1}^\Delta(X)$ by definition.  Hence we have  $\bigl(E^{{\rm BM}, 1}_{\ast,0}(X), d^1_{\ast, 0}\bigr)\cong\bigl(\Omega_\ast^{\rm BM}(X),\del^{\rm path}_\ast\bigr)$ and obtain  $E^{{\rm BM}, 2}_{p, 0}(X) \cong  {\rm PH}^{\rm BM}_p(X)$. The uniformly finite case is analogous. This completes the proof.

\end{proof}

\subsubsection{Some consequences}\label{consequences}

The following shows the small-scale to large-scale interpolation of  the magnitude--path spectral sequence. It shows that an isomorphism on a finite page of the spectral sequence is sufficient to detect an isomorphism on coarse homology.
\begin{prp}
Let $X$ and $Y$ be locally finite graphs. If a filtered chain map $f : C_\ast(X) \to C_\ast(Y)$ induces isomorphisms on $E^{{\rm BM}, r}_{p, q}$-terms for some $r\geq 0$ and for all $p, q$, then it induces an isomorphism on the coarse homology $f_\ast : H_\ast(X) \cong H_\ast(Y)$. Likewise, when $X$ and $Y$ have bounded geometry, if a filtered chain map $f : C^{\rm uf}_\ast(X) \to C^{\rm uf}_\ast(Y)$ induces isomorphisms on $E^{\mathrm{uf},r}_{p,q}$-terms for some $r\geq 0$ and for all $p, q$, then it induces an isomorphism on uniformly finite homology.
\end{prp}
\begin{proof}
This is an immediate consequence of Proposition  \ref{comparison}, since the length filtrations are bounded below and exhaustive.
\end{proof}

Next, although the terms $E^r_{p,q}$ need not stabilize at any finite page in general, the following result shows that $E^r_{0,0}$ does stabilize at the second page.

\begin{prp}\label{h0}
For a locally finite graph $X$, we have
\begin{align*}
H_0(X)\cong E^{{\rm BM}, \infty}_{0,0}(X)&\cong E^{{\rm BM}, r}_{0,0}(X) \cong E^{{\rm BM}, 2}_{0,0}(X) \cong {\rm PH}^{\rm BM}_0(X), \\
H_0^{\rm uf}(X)\cong E^{{\rm uf}, \infty}_{0,0}(X)&\cong E^{{\rm uf}, r}_{0,0}(X) \cong E^{{\rm uf}, 2}_{0,0}(X) \cong {\rm PH}^{\rm uf}_0(X), 
\end{align*}
for all $r \geq 2$, where the latter is defined when $X$ has bounded geometry.
\end{prp}
\begin{proof}
We give a proof for $E^{{\rm BM}, r}_{p,q}(X)$, where we denote it by $E^r_{p, q}$ in the following. The proof for uniformly finite version is analogous. From the construction, we have $Z^r_{0, 0} = F_0C_0(X) = C_0(X)$ for all $r\geq1$, and hence $d^r_{0, 0} = 0 : E^r_{0,0}\to  E^r_{-r, r-1}$ for $r\geq1$. Also, as shown in \cite[Corollary~4.5]{LS}, the degree $1$ magnitude homology ${\rm MH}_{1, \ell}(X)$ is freely generated by ordered adjacent pairs of vertices, in particular, we have ${\rm MH}^x_{1, \ell}(X) = 0$ for all $x \in X$ and $\ell\geq 2$. This implies that 
\begin{align*}
E^1_{r, -r+1} &= {\rm MH}^{\rm BM}_{1, r}(X) = 0,
\end{align*}
for all $r\geq2$ by Proposition \ref{prod}. Hence we obtain $E^r_{r, -r+1} = 0$ and $d^r_{r, -r+1}= 0 : E^r_{r, -r+1} \to E^r_{0, 0}$ for all $r\geq2$. Thus $E^r_{0, 0}$ stabilizes at $r=2$, which implies $B^{r-1}_{0, 0} = B^1_{0, 0}$ for $r\geq 2$ since $E^r_{0, 0} = Z^r_{0, 0}/B^{r-1}_{0, 0} = C_0(X)/B^{r-1}_{0, 0}$.  Thus we obtain $E^2_{0, 0} = E^\infty_{0, 0}$, and the isomorphism with the path homology is given by Proposition \ref{pathhomology}. Also, we have $E^{\infty}_{0,0} \cong F_0H_0/F_{-1}H_0=H_0(X)$. This completes the proof.
\end{proof}

\begin{rem}
Proposition \ref{h0} is closely related to the characterization of amenability by uniformly finite homology. Let $\Gamma$ be a finitely generated group and let $X={\rm Cay}(\Gamma,S)$ be a Cayley graph with respect to a finite generating set $S$. As in \cite[Definition~7.2.3]{NY}, the constant function $1\in C^{\rm uf}_0(X)$ determines the fundamental class $[X]\in H^{\mathrm{uf}}_0(X)$. By \cite[Theorem~7.2.4, Proposition~7.3.4]{NY}, we have
\[
H^{\rm uf}_0(X)=0 \Leftrightarrow [X]=0\Leftrightarrow \Gamma\text{ is non-amenable}.
\]
Hence Proposition \ref{h0} implies that
\[
{\rm PH}^{\rm uf}_0(X)=0 \Leftrightarrow E^{\mathrm{uf},r}_{0, 0}(X)=0\Leftrightarrow \Gamma\text{ is non-amenable},
\]
for all $r\geq 2$. Thus amenability of a finitely generated group is already detected on the second page of the uniformly finite magnitude--path spectral sequence.
\end{rem}

Now we give a sufficient condition for the magnitude--path spectral sequences to degenerate at the $E^2$-page.

\begin{prp}\label{e2deg}
Let $X$ be a locally finite graph which is diagonal in the sense of ordinary magnitude homology. Then the magnitude--path spectral sequences $E^{{\rm BM}, \ast}_{\ast, \ast}$ and $E^{{\rm uf}, \ast}_{\ast, \ast}$ degenerate at $E^2$-page, and we have
\begin{align*}
H_p(X)&\cong {\rm PH}^{\rm BM}_{p}(X), \\
H^{\rm uf}_p(X)&\cong {\rm PH}^{\rm uf}_{p}(X),
\end{align*}
for all $p\geq0$, where the latter is defined when $X$ has bounded geometry.
\end{prp}
\begin{proof}
By Proposition \ref{diagonality}, the diagonality of the ordinary magnitude homology ${\rm MH}_{\ast,\ast}(X)$ implies the diagonality of ${\rm MH}^{\rm BM}_{\ast,\ast}(X)$ and, when $X$ has bounded geometry, of ${\rm MH}^{\rm uf}_{\ast,\ast}(X)$. Hence Proposition \ref{weaklyconverge} gives
\[
E^{{\rm BM},1}_{p,q}(X)
\cong {\rm MH}^{\rm BM}_{p+q,p}(X)=0
\]
for $q\neq 0$, and similarly for the uniformly finite version. Thus the $E^1$-page is concentrated on the $p$-axis. Consequently, the $E^2$-page is also concentrated on the $p$-axis, and all differentials $d^r$ for $r\geq 2$ vanish. Therefore the spectral sequences degenerate at the $E^2$-page. Since we have $Z^r_{p,0}=F_pC_p^\Delta(X)\cap\ker\partial_p^\Delta$ and $Z^{r-1}_{p-1,1}=0$ for $r\geq2$, the same argument as in Proposition \ref{h0} applies, and we obtain $E^2_{p, 0}(X) \cong E^\infty_{p, 0}(X)$ for all $p\geq0$. Since the spectral sequences  converge and, in each total degree $p$, the only nonzero $E^\infty$-term is $E^\infty_{p,0}$, there is no extension problem. Hence the desired isomorphism is obtained by Proposition \ref{pathhomology}. This completes the proof.
\end{proof}

We give examples of diagonal graphs to which Proposition \ref{e2deg} can be applied. Recall that a connected graph $X$ is called {\it quasi-median} if it is weakly modular and contains neither $K_{1, 1, 2}$ nor $K_{2, 3}$ as an induced subgraph. See \cite{Genevois} for details. The following is proved analogously to \cite[Proposition~2.2]{BK}.

\begin{lem}
Every quasi-median graph is diagonal.
\end{lem}
\begin{proof}
By \cite{Chastand}, every quasi-median graph is a retract of a Hamming graph, namely a Cartesian product of complete graphs. Note that every complete graph is diagonal, and diagonality is preserved under Cartesian products by \cite[Example~2.5, Proposition 7.3]{HW}. By \cite[Propositions~2.5, 2.6]{BK}, the diagonality is also preserved under retracts, and filtrations by isometric subgraphs. Therefore every Hamming graph is diagonal, and so is every quasi-median graph.
\end{proof}

Genevois \cite{Genevois} gives an introduction to quasi-median graphs and describes several examples,
including graphs associated with graph products, wreath products, and diagram products. In particular, the Cayley graphs of {\it graph products}, including {\it right-angled Artin groups and right-angled Coxeter groups} are quasi-median \cite[Theorem~3.2]{Genevois}. The graph product is constructed as follows: let $\Gamma$ be a graph and let $\mathcal{G}=\{G_v\}_{v\in V(\Gamma)}$ be a collection of groups. The graph product $\Gamma\mathcal{G}$ is a group defined as 
\[
\langle G_v, v \in V(\Gamma) \mid [G_u, G_v] = 1 \text{ for } (u, v) \in E(\Gamma) \rangle.
\]
Let ${\rm Cay}(\Gamma,\mathcal{G})$ be the Cayley graph of $\Gamma\mathcal{G}$ associated to the generating set $\cup_{v \in V(\Gamma)}G_v\setminus \{1\}$. When it is locally finite, Proposition
\ref{e2deg} applies and gives
\begin{align*}
H_\ast({\rm Cay}(\Gamma,\mathcal{G})) \cong {\rm PH}^{\rm BM}_\ast({\rm Cay}(\Gamma,\mathcal{G})), \\
H_\ast^{\rm uf}({\rm Cay}(\Gamma,\mathcal G))\cong {\rm PH}^{\rm uf}_*({\rm Cay}(\Gamma,\mathcal G)),
\end{align*}
where the latter is defined when the Cayley graph has bounded geometry. Note that these Cayley graphs are not necessarily locally finite. In particular, the Cayley graphs of right-angled Artin groups with the generating sets considered here are not locally finite. However, it is well known that their Cayley graphs with respect to the standard generating sets are median graphs. Hence the above isomorphisms also hold for these Cayley graphs when the defining graphs are finite.

In contrast to the above examples, the following shows examples such that the spectral sequence needs not degenerate but the path homologies are very simple. Recall that {\it the girth} of a graph $X$ is the least length of cycles in $X$. The following is essentially proved in \cite[Proposition~8.13]{A}.
\begin{prp}\label{girth}
Let $X$ be a locally finite graph of girth at least $5$. Then we have 
\begin{align*}
{\rm PH}^{\rm BM}_n(X)&=0, \\
{\rm PH}^{\rm uf}_n(X) & = 0,
\end{align*}
for all $n\geq 2$, where the latter is defined when $X$ has bounded geometry. 
\end{prp}

\begin{proof}
By the proof of \cite[Proposition 8.13]{A}, the complex
\[
\cdots \to {\rm MH}_{n,n}(X) \xrightarrow{d^1} {\rm MH}_{n-1,n-1}(X) \to \cdots
\]
is exact for $n\geq2$. Also, by \cite[Corollary~1.2]{AHK}, after truncating in degrees below $1$, this complex decomposes as a direct sum of complexes 
\[
\cdots \to {\rm MH}^{\{x, y\}}_{n,n}(X) \xrightarrow{d^1} {\rm MH}^{\{x, y\}}_{n-1,n-1}(X) \to \cdots,
\]
where ${\rm MH}^{\{x, y\}}_{n,n}(X) = \Z\langle (x, y, x, \dots), (y, x , y, \dots)\rangle$. Taking the product of its components and using Proposition \ref{prod}, we obtain the exactness of
\[
\cdots \to {\rm MH}^{\rm BM}_{n,n}(X)\xrightarrow{d^1}{\rm MH}^{\rm BM}_{n-1,n-1}(X)\to \cdots
\]
 for $n\geq2$.  Hence Proposition~\ref{pathhomology} implies ${\rm PH}^{\rm BM}_n(X)=0$ for $n\geq2$. For the uniformly finite version, note that the above subcomplexes are all norm-preservingly isomorphic, with respect to the $\ell^\infty$-norms induced by their natural bases. Since each differential is strict by Lemma \ref{freeabel},  the corresponding family of differentials is uniformly strict by Lemma \ref{fag}. Thus, by Lemma \ref{prodchain}, we have that the restricted product of these complexes is exact in degrees at least $2$. By Proposition \ref{prod} and Proposition \ref{pathhomology}, we obtain ${\rm PH}^{\rm uf}_n(X)=0$ for $n\geq2$. This completes the proof.
\end{proof}

\subsection{Functoriality and invariance}\label{sec:5.3}

\subsubsection{Functoriality}
\begin{df}
Let $f : X \to Y$ be a map between locally finite graphs.
\begin{enumerate}
\item We say $f$ is {\it a coarse map} if, for all $R>0$, there exists $S>0$ such that $d_X(x, x')\leq R$ implies $d_Y(f(x), f(x'))\leq S$ for all $x, x' \in X$.
\item We say $f$ is {\it proper} if every fiber of $f$ is a finite set.
\item We say $f$ has {\it uniformly finite fibers} if it is proper and the cardinalities of the fibers are uniformly bounded.
\end{enumerate}
\end{df}

\begin{rem}\label{coarselipschitz}
For graphs equipped with the path metric, a map $f : X \to Y$ is a coarse map if and only if we have a constant $A>0$ such that $d(f(x), f(x')) \leq Ad(x, x')$ for all $x, x' \in X$. In fact, if $d(x, x') = 1$ implies $d(f(x), f(x')) \leq A$ for all $x, x' \in X$, then we have $d(f(x), f(x')) \leq Ad(x, x')$ for all $x, x' \in X$. Note that a coarse map $f$ is 1-Lipschitz if and only if we can take $A = 1$.
\end{rem}
\begin{prp}\label{functorial}
Let $f : X \to Y$ be a map between locally finite graphs.
\begin{enumerate}
\item If $f$ is a proper coarse map, it induces a morphism of simplicial abelian groups $f_\ast : C_\bullet (X) \to C_\bullet(Y)$. Moreover, if $f$ is also 1-Lipschitz, $f_\ast$ preserves the filtration $F_\bullet$, hence it induces a morphism of spectral sequences $f^r_{p, q} : E^{{\rm BM}, r}_{p, q}(X) \to E^{{\rm BM}, r}_{p, q}(Y)$ for all $p, q, r$.
\item Suppose that $X$ and $Y$ have bounded geometry. If $f$ is a coarse map with uniformly finite fibers, it induces a morphism of simplicial abelian groups $f_\ast : C^{\rm uf}_\bullet (X) \to C^{\rm uf}_\bullet(Y)$. Moreover, if $f$ is also 1-Lipschitz, $f_\ast$ preserves the filtration $F_\bullet$, hence it induces a morphism of spectral sequences $f^r_{p, q} : E^{{\rm uf}, r}_{p, q}(X) \to E^{{\rm uf}, r}_{p, q}(Y)$ for all $p, q, r$.
\end{enumerate}
\end{prp}
\begin{proof}
\begin{enumerate}
\item For $\varphi \in C_n(X)$, we define $f_\ast\varphi : Y^{n+1} \to \Z$ by
\[
f_\ast\varphi (y_0,\dots, y_n) = \sum_{(f(x_0), \dots, f(x_n)) = (y_0,\dots, y_n)} \varphi(x_0, \dots, x_n),
\]
which is well-defined by the properness assumption on $f$. Since $\varphi$ has controlled support, there exists a constant $R>0$ such that $d(x_i, x_j) \leq R$ for all $i, j$ if $(x_0, \dots, x_n) \in {\rm Supp}\varphi$. Then, since $f$ is a coarse map, there exists a constant $S>0$ such that $d(y_i, y_j) \leq S$ for all $i, j$ if $(y_0, \dots, y_n) \in {\rm Supp}f_\ast \varphi$, which implies that $f_\ast \varphi$ has controlled support, and hence $f_\ast \varphi \in C_n(Y)$. The compatibility with the face and the degeneracy maps is straightforward. Also, it is immediate to verify that $f_\ast$ preserves the filtration $F_\bullet$ if it is 1-Lipschitz.
\item For $\varphi \in C^{\rm uf}_n(X)$, we define $f_\ast\varphi : Y^{n+1} \to \Z$ by the same formula as in (1). Then it is well-defined and has controlled support by the same argument. Further, since $f$ has uniformly finite fibers, there exists a constant $N>0$ such that $\# f^{-1}(y) \leq N$ for all $y \in Y$, which implies 
\[
|f_\ast\varphi (y_0,\dots, y_n)| \leq N^{n+1}|\varphi|,
\]
hence $f_\ast \varphi$ is bounded and $f_\ast \varphi \in C^{\rm uf}_n(Y)$. The remaining argument is the same as in (1). This completes the proof.
\end{enumerate}
\end{proof}

\subsubsection{Invariance}
\begin{df}\label{rhomotopy}
Let $f, g : X \to Y$ be proper 1-Lipschitz maps between locally finite graphs. 
\begin{enumerate}
\item For $r\geq0$, we say $f$ and $g$ are {\it $r$-close} if we have $d_Y(f(x), g(x)) \leq r$ for all $x \in X$. 
\item We denote the equivalence relation generated by the $r$-closeness by $\sim_r$, and we say {\it $f$ is $r$-homotopic to $g$} if $f\sim_r g$.
\item We say $X$ and $Y$ are {\it $r$-homotopy equivalent} if there exist proper 1-Lipschitz maps $f : X \to Y$ and $g : Y \to X$ such that $g\circ f \sim_r {\rm id}_X$ and $f\circ g \sim_r {\rm id}_Y$. 
\end{enumerate}
\end{df}

\begin{prp}\label{invariance}
For a locally finite graph $X$, $E^{{\rm BM}, r+1}_{\bullet, \bullet}(X)$ is an $r$-homotopy invariant. Namely, if we have proper 1-Lipschitz maps $f, g : X \to Y$ between locally finite graphs, we have
\[
f^{r+1}_{p, q} = g^{r+1}_{p, q} : E^{{\rm BM}, r+1}_{p, q}(X) \to E^{{\rm BM}, r+1}_{p, q}(Y),
\]
for all $p, q$ provided that $f\sim_r g$. In particular, if $X$ and $Y$ are $r$-homotopy equivalent, then we have $ E^{{\rm BM}, R}_{p, q}(X) \cong E^{{\rm BM}, R}_{p, q}(Y)$ for all $R>r$ and all $p, q$. Also, if $X$ and $Y$ have bounded geometry, we have the analogous statement for $E^{{\rm uf}, r+1}_{\bullet, \bullet}(X)$ with an additional condition that $f, g$ have uniformly finite fibers.
\end{prp}

The proof of Proposition \ref{invariance} is essentially the same as those of  \cite[Proposition~3.9]{Cir} and \cite[Theorem~4.37]{Asa}, both of which refer to the classical argument in \cite[p.321]{CE}. Here we record a version adapted to our notation.

\begin{lem}[\cite{CE}p. 321]\label{CE}
Let $C_\ast$ and $D_\ast$ be chain complexes with increasing filtrations $F_\bullet C_\ast$ and $F_\bullet D_\ast$. Let $f,g : C_\ast \to D_\ast$ be filtered chain maps, and suppose that there exists a chain homotopy $h : C_\ast \to D_{\ast+1}$ from $f$ to $g$. Assume moreover that
\[
h(F_pC_n)\subset F_{p+r}D_{n+1}
\]
for all $p$ and $n$. Then we have $E^{r+1}(f)=E^{r+1}(g) : E^{r+1}(C) \to E^{r+1}(D)$. 
\end{lem}
\begin{proof}
Let 
\[
x \in Z^{r+1}_{p,q}(C)
=
F_pC_{p+q}\cap
\del_{p+q}^{-1}\bigl(F_{p-r-1}C_{p+q-1}\bigr).
\]
Note that, we have 
\begin{align*}
\del h(x)&\in \del(F_{p+r}D_{p+q+1}), \\
 h\del(x)&\in F_{p-1}D_{p+q},
\end{align*}
from $h(F_pC_{p+q})\subset F_{p+r}D_{p+q+1}$ and $\del (x) \in F_{p-r-1}C_{p+q-1}$. Since we have $\del h(x) = (f-g)(x) - h\del(x)$ and $f-g$ preserves filtrations, we obtain $\del h(x)\in F_pD_{p+q}$, and hence 
\[
\del h(x)\in F_pD_{p+q}\cap \del\bigl(F_{p+r}D_{p+q+1}\bigr)
=
B^r_{p,q}(D).
\]
Also, from $(f-g)\del (x)\in F_{p-r-1}D_{p+q-1}$ and that we have $\del h\del (x) = (f-g)\del (x)$ by $\del^2=0$, we obtain  
\[
h\del(x)\in
F_{p-1}D_{p+q}
\cap
\del^{-1}\bigl(F_{p-r-1}D_{p+q-1}\bigr) = Z^{r}_{p-1, q+1}.
\]
Therefore, we obtain 
\[
(f-g)(x) = 0 \in E^{r+1}_{p, q} = Z^{r+1}_{p,q}/\bigl(Z^r_{p-1, q+1} + B^r_{p, q}\bigr).
\]
This completes the proof.
\end{proof}

\begin{proof}[Proof of Proposition \ref{invariance}]
Since $r$-homotopy is generated by $r$-closeness, it suffices to prove the assertion when $f$ and $g$ are $r$-close. We show that $f^{r+1}_{p, q} = g^{r+1}_{p, q}$ if the proper 1-Lipschitz maps $f, g : X \to Y$ are $r$-close. For  $\varphi \in C_n(X)$, we define a map $h_n\varphi : Y^{n+2} \to \Z$ by
\[
h_n\varphi(y_0, \dots, y_{n+1}) = \sum_{i=0}^n(-1)^i\sum_{\substack{(y_0,\dots, y_i) = (f(x_0), \dots, f(x_i))\\ (y_{i+1},\dots, y_{n+1}) = (g(x_i), \dots, g(x_n))}} \varphi(x_0,\dots, x_n),
\]
which is well-defined by the properness assumption on $f$ and $g$.  Since $\varphi$ has controlled support, there exists a constant $\ell>0$ such that $\varphi \in F_\ell C_n(X)$, namely $L(x_0, \dots, x_n) \leq \ell$ for $(x_0, \dots, x_n) \in {\rm Supp} \varphi$. Then, since $f$ and $g$ are 1-Lipschitz and $r$-close, we have 
\[
L(y_0, \dots, y_{n+1}) = \sum_{k=0}^{i-1} d(f(x_{k}), f(x_{k+1})) + d(f(x_i), g(x_i)) +  \sum_{k=i}^{n-1} d(g(x_{k}), g(x_{k+1})) \leq \ell+r,
\]
for $(y_0, \dots, y_{n+1}) \in {\rm Supp}h_n\varphi$, which implies that $h_n\varphi \in F_{\ell+r}C_{n+1}(Y)$, in particular that $h_n\varphi$ has controlled support. Also, it is straightforward to verify that
\[
\del_{n+1}h_n\varphi + h_{n-1}\del_n\varphi = g_\ast\varphi - f_\ast\varphi .
\]
Therefore, the correspondences $\{h_n\}_n$ define a chain homotopy $h_\ast : C_\ast(X) \to C_{\ast+1}(Y)$ with 
\[
h_\ast(F_\ell C_\ast(X)) \subset F_{\ell+r}C_{\ast+1}(Y).
\]
Then Lemma \ref{CE} implies the statement. For the case of uniformly finite version, we have $|h_n\varphi|
\leq (n+1)N^{n+1}|\varphi|$, which implies that  $h_n\varphi$ is bounded, provided that $X$ and $Y$ have bounded geometry, $f$ and $g$ have uniformly finite fibers and $\varphi$ is bounded.  The remaining argument is the same. Note that, under the bounded geometry assumption, any map $r$-close to a map with uniformly finite fibers also has uniformly finite fibers. This completes the proof.
\end{proof}

\subsubsection{Relation to quasi-isometries}

\begin{df}\label{coarsemap}
Let $X$ and $Y$ be locally finite graphs. We say a map $f : X \to Y$ is {\it a quasi-isometry} if there exist constants $A, B, C > 0$ satisfying the following:
\begin{itemize}
\item For all $x, x' \in X$,
\[
\frac{1}{A}d_X(x, x') - B\leq d_Y(f(x), f(x')) \leq Ad_X(x, x') + B.
\]
\item For all $y \in Y$, there is a $x \in X$ such that $d_Y(f(x), y) \leq C$.
\end{itemize}
\end{df}

\begin{rem}\label{qisoproper}
Every quasi-isometry $f : X \to Y$ is proper. In fact, for $x, x' \in X$ with $f(x) = f(x')$, we have $d_X(x, x') \leq AB$, namely every fiber of $f$ has finite diameter. Since bounded subsets of a locally finite graph are finite, every fiber of $f$ is finite, and hence $f$ is proper.  Moreover, if $X$ also has bounded geometry, quasi-isometry implies  the uniform boundedness of fibers of $f$.
\end{rem}

\begin{rem}\label{coarseqiso}
We say metric spaces $X$ and $Y$ are {\it coarsely equivalent} if there are proper coarse maps $f : X \to Y$ and $g : Y \to X$ such that $d(g\circ f(x), x) < C$ and $d(f\circ g(y), y) < C$ for some constant $C>0$. We say $f$ and $g$ are {\it coarse equivalence maps}. As long as we consider locally finite graphs, a map is a coarse equivalence if and only if it is a quasi-isometry (\cite[Corollary~1.4.14]{NY}).
\end{rem}

Note that proper coarse maps $f, g : X \to Y$ induce chain maps $f_\ast, g_\ast : C_\ast(X) \to C_\ast(Y)$. If $f$ and $g$ are $r$-close for some $r$, it is straightforward to construct a chain homotopy between $f_\ast$ and $g_\ast$. Hence Remarks  \ref{qisoproper} and \ref{coarseqiso} imply the following (see also \cite[Theorem~7.15]{NY}).

\begin{prp}\label{coarseinvariance}
For a locally finite graph $X$, the coarse homology $H_\ast(X)$ is a coarse invariant. Also, if $X$ has bounded geometry, the uniformly finite homology $H_\ast^{\rm uf}(X)$ is a quasi-isometry invariant within the class of graphs of bounded geometry.
\end{prp}

Note that a homomorphism $f_\ast : C_\ast(X) \to C_\ast(Y)$ induced by a proper coarse map $f : X \to Y$ does not preserve the filtration $F_\bullet$ in general. Thus coarse or quasi-isometric maps do not in general induce morphisms of the spectral sequences preserving the bidegree,  while the following shows that it is invariant under special type of quasi-isometries, namely an $r$-homotopy equivalence.
\begin{lem}
Let $X$ and $Y$ be locally finite graphs. Then an $r$-homotopy equivalence $f : X \to Y$ is a quasi-isometry.
\end{lem}
\begin{proof}
Since $f$ is 1-Lipschitz, we have $d(f(x), f(x')) \leq d(x, x')$ for all $x, x' \in X$. Let $g : Y \to X$ be an $r$-homotopy inverse. By the definition of $r$-homotopy, there exists $N\geq 0$ such that $d(gf(x), x)\leq Nr$ for all $x \in X$, and $g$ is a 1-Lipschitz map. Hence we have
\begin{align*}
d(x, x') &\leq d(x, gf(x)) + d(gf(x), gf(x')) + d(gf(x'), x') \\
&\leq 2Nr + d(f(x), f(x')),
\end{align*}
which implies that $d(x, x') - 2Nr \leq  d(f(x), f(x')) \leq d(x, x')$. Also, by the definition of $r$-homotopy, there exists $M\ge0$ such that $d(y, fg(y))\leq Mr$ for all $y \in Y$. Therefore $f$ is a quasi-isometry.  This completes the proof.
\end{proof}

\subsection{An example of non-stabilization and a uniformly finite coarse class}\label{sec:5.4}

The purpose of this subsection is to illustrate the distinction between  convergence and finite-page stabilization discussed in Remark \ref{stability}, and at the same time to exhibit a coarse-geometric phenomenon detected by uniformly finite homology. We construct a graph of bounded geometry containing a sequence of finite cycles such that each individual cycle becomes trivial at a sufficiently large scale, while the scales required to kill these cycles tend to infinity. Consequently, their pointwise limit cannot be filled with uniformly bounded control and represents a non-trivial class in uniformly finite homology. From the viewpoint of the magnitude--path spectral sequence, this lack of a uniform filling scale appears as the failure of $E^{{\rm uf}, r}_{1, 0}$  to stabilize at any finite page.

Let us define a graph $X_n$ inductively on $n$ as follows. Let $C_{2n+1}$ be a $(2n+1)$-cycle graph with vertices $\{v_i^{2n+1}\mid i\in\Z/(2n+1)\}$ indexed clockwise. We define $X_1 = C_3$ and $X_{n+1} = X_n \vee_{0} e_n \vee_{1} C_{2n+3}$, where $e_n$ denotes an edge graph with vertices $v^n_0$ and $v^n_1$, and $\vee_i$ for $i=0, 1$ denotes the wedge sum identifying vertices $v^{2(n+i)+1}_i$ and $v^n_i$. Obviously we have an inclusion $X_m \to X_n$ for $m\leq n$, and we define an infinite graph $X = {\rm colim}_nX_n$. The graph $X$ has uniformly bounded degree, and hence has bounded geometry. Now we consider functions $f_n : X^2 \to \Z$ defined by 
\[
f_n(x, y) = \begin{cases} 1 & (x,y) = (v^{2k+1}_i, v^{2k+1}_{i+1}) \text{ for some } k\leq n \text{ and } i,\\ 0 & \text{otherwise}, \end{cases}
\]
and their pointwise limit $f = \lim_n f_n$. Then it is easy to verify that $f_n, f \in C^{\rm uf}_1(X)$ and $\del_1 f_n= \del_1 f =  0$, hence we have $[f_n], [f] \in E^{{\rm uf}, r}_{1, 0}(X)$. 

Note that $X_n$ is $n$-homotopy equivalent to a point. Indeed, we may first collapse the terminal cycle $C_{2n+1}$ to its attaching vertex; this map is $n$-close to the identity since ${\rm diam} C_{2n+1}=n$. We may then collapse the adjacent edge, which is $1$-close to the identity. Thus $X_n$ is $n$-homotopy equivalent to $X_{n-1}$, and induction together with the  $1$-homotopy equivalence  \(X_1=C_3\simeq \ast$ gives an $n$-homotopy equivalence between $X_n$ and a point. 
Hence Proposition \ref{invariance} gives $[f_n]=0\in E^{{\rm uf},r}_{1, 0}(X_n)$ for $r>n$. Applying the homomorphism induced by the inclusion $X_n\hookrightarrow X$, we obtain the same vanishing in $E^{{\rm uf},r}_{1, 0}(X)$. 

On the other hand, Proposition \ref{EMV} below, applied to the decomposition of $X$ at the bridge following $X_n$, shows that $E^{{\rm uf}, r}_{1, 0}(X_n)$ is a direct summand of $E^{{\rm uf}, r}_{1, 0}(X)$. Applying Proposition \ref{EMV} also to the decomposition of $X_n$ at the bridge preceding $C_{2n+1}$, we see that $E^{{\rm uf}, r}_{1, 0}(C_{2n+1})$ is a direct summand of $E^{{\rm uf}, r}_{1, 0}(X_n)$. The $C_{2n+1}$-component of $[f_n]$ is nonzero in $E^{{\rm uf}, r}_{1, 0}(C_{2n+1})$ for $r\leq n$ by Lemma \ref{NT} below. Hence we obtain
\[
[f_n]\ne0\in E^{{\rm uf},r}_{1, 0}(X),
\]
for $r\leq n$. Therefore, for every $r\geq1$, the class $[f_r]$ is nonzero in $E^{{\rm uf}, r}_{1, 0}(X)$ but vanishes in $E^{{\rm uf}, r+1}_{1, 0}(X)$. Hence $E^{{\rm uf}, r}_{1, 0}(X)$ does not stabilize at any finite page.

Suppose that we have $[f] = 0$ in $H^{\rm uf}_1(X)$. Then there is $b \in F_RC_2^{\rm uf}(X)$ with $\del_2b = f$ for some $R>0$, hence we have $[f] = 0$ in $E^{{\rm uf}, r}_{1, 0}(X)$ for sufficiently large $r$. However, by Proposition \ref{EMV}, $E^{{\rm uf}, r}_{1, 0}(X)$ has $E^{{\rm uf}, r}_{1, 0}(X_n)$ as a direct summand for all $n$, and we have a decomposition $[f] = [f_n] \oplus [g]$. Since the $C_{2n+1}$-component of $[f_n]$ in $E^{{\rm uf}, r}_{1, 0}(X)$ is the non-zero class of Lemma \ref{NT} for $n\geq r$, this is a contradiction. Therefore we obtain a non-trivial uniformly finite homology class $[f] \in H^{\rm uf}_1(X)$. This example demonstrates the case mentioned in Remark \ref{stability}, namely $E^{{\rm uf}, r}_{1, 0}(X)$ does not stabilize while the spectral sequence  converges.
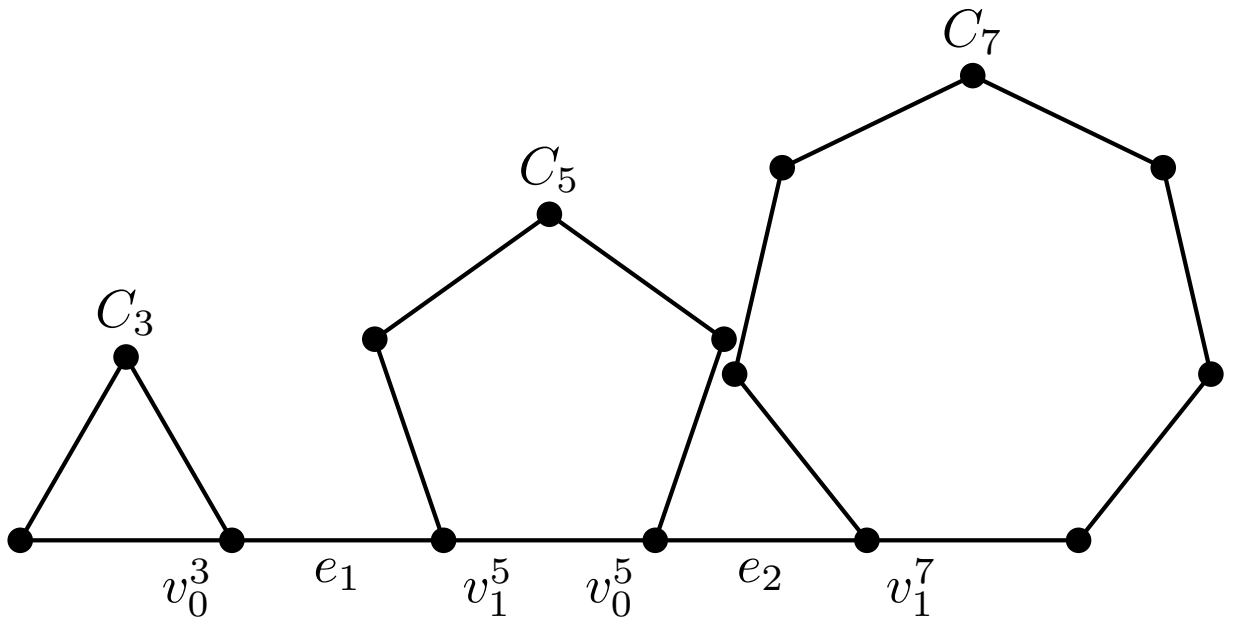
\begin{figure}[H]
\center
\begin{tikzpicture}[
    scale=0.7,
    vertex/.style={circle, fill=black, inner sep=1.7pt},
    edge/.style={thick}
]

\coordinate (A0) at (-1,0);
\coordinate (A1) at (1,0);
\coordinate (A2) at (0,1.732);

\draw[edge] (A0)--(A1)--(A2)--cycle;

\node[vertex] at (A0) {};
\node[vertex] at (A1) {};
\node[vertex] at (A2) {};

\node at (0,2.15) {$C_3$};


\coordinate (B0) at (3,0);
\coordinate (B1) at (5,0);

\coordinate (B2) at (5.65,1.90);
\coordinate (B3) at (4,3.08);
\coordinate (B4) at (2.35,1.90);

\draw[edge]
    (B0)--(B1)--(B2)--(B3)--(B4)--cycle;

\node[vertex] at (B0) {};
\node[vertex] at (B1) {};
\node[vertex] at (B2) {};
\node[vertex] at (B3) {};
\node[vertex] at (B4) {};

\node at (4,3.5) {$C_5$};


\coordinate (C0) at (7,0);
\coordinate (C1) at (9,0);

\coordinate (C2) at (10.25,1.57);
\coordinate (C3) at (9.80,3.52);
\coordinate (C4) at (8,4.39);
\coordinate (C5) at (6.20,3.52);
\coordinate (C6) at (5.75,1.57);

\draw[edge]
    (C0)--(C1)--(C2)--(C3)--(C4)--(C5)--(C6)--cycle;

\node[vertex] at (C0) {};
\node[vertex] at (C1) {};
\node[vertex] at (C2) {};
\node[vertex] at (C3) {};
\node[vertex] at (C4) {};
\node[vertex] at (C5) {};
\node[vertex] at (C6) {};

\node at (8,4.8) {$C_7$};

\draw[edge] (A1)--(B0);
\draw[edge] (B1)--(C0);

\node[below] at (2,0) {$e_1$};
\node[below] at (6,0) {$e_2$};

\node[below left]  at (A1) {$v_0^3$};
\node[below right] at (B0) {$v_1^5$};

\node[below left]  at (B1) {$v_0^5$};
\node[below right] at (C0) {$v_1^7$};

\end{tikzpicture}
\caption{The picture of $X_3$.}
\end{figure}

\begin{lem}\label{NT}
Let $\varphi  \in C_1(C_{2n+1})$ defined by 
\[
\varphi (x, y) = \begin{cases} 1 & (x, y) = (v^{2n+1}_i, v^{2n+1}_{i+1}) \text{ for some } i, \\ 0 & \text{otherwise}.\end{cases}
\]
Then we have $[\varphi]\neq 0 \in E^{r}_{1, 0}(C_{2n+1})$ for $r\leq n$. Here, since $C_{2n+1}$ is finite, we omit the superscript ${\rm BM}$ or ${\rm uf}$ from its magnitude--path spectral sequence.
\end{lem}
\begin{proof}
Since a class which vanishes on the $r$-th page also vanishes on every subsequent page, it suffices to prove the assertion for $r=n$. Suppose that we have $[\varphi]= 0 \in E^{n}_{1, 0}(C_{2n+1})$. Then there exist $z \in F_0C_1(C_{2n+1})$ and $b \in F_nC_2(C_{2n+1})$ such that $z + \del_2b = \varphi$. Now we define a map $w : V(C_{2n+1})^2 \to \Z$ by $w(x, y) = d(x, y)$, and we also define a map $W : C_1(C_{2n+1}) \to \Z$ by 
\[
W(\psi) = \sum_{(x, y) \in V(C_{2n+1})^2} \psi(x, y)\cdot w(x, y).
\]
Note that, if $L(x, y, z) \leq n$, then we have either $d(x, y) + d(y, z) = d(x, z), d(x, z) + d(y, z) = d(x, y)$ or $d(x, y) + d(x, z) = d(y, z)$. In each case, we have 
\[
w(y, z) - w(x, z) + w(x, y) \in 2\Z,
\]
which implies
\begin{align*}
W(\del_2 b) &= \sum_{(x, y) \in V(C_{2n+1})^2} \del_2b(x, y)\cdot w(x, y) \\
 &= \sum_{(x, y, z) \in {\rm Supp}b} b(x, y, z)\cdot (w(y, z) - w(x, z) + w(x, y)) \\
 &\in 2\Z,
\end{align*}
since we have $b \in F_nC_2(C_{2n+1})$. Also, we have $W(z) = 0$ since we have $z \in F_0C_1(C_{2n+1})$. On the other hand, we have $W(\varphi) = 2n+1 \not\in 2\Z$, which is a contradiction. This completes the proof.
\end{proof}

\section{Mayer--Vietoris  sequences}\label{MVS}

In this section, we establish Mayer--Vietoris sequences associated with projecting decompositions. We begin in Section \ref{sec6.1} by extending the split Mayer--Vietoris sequence of Hepworth--Willerton for ordinary magnitude homology to the Borel--Moore and uniformly finite versions. In Section \ref{sec6.2} , the resulting short exact sequence on the $E^1$-page gives Mayer--Vietoris long exact sequences for path homology. Under suitable finiteness assumptions on the fibers of the projecting map, the splitting is compatible with the spectral sequence differentials and yields split short exact sequences on every page, as shown in Section \ref{sec6.3} . Finally, in Section \ref{sec6.4} , we use a mapping cone construction together with the comparison theorem to obtain Mayer--Vietoris long exact sequences for coarse and uniformly finite homology without any finiteness assumption on the fibers of the projecting map.

\subsection{Mayer--Vietoris sequence for magnitude homology}\label{sec6.1}

The following proposition is based on Hepworth--Willerton \cite[Theorem~6.5, 6.6]{HW}, but we refine it to components with a fixed initial point. Further, although their convention assumes finiteness of graphs, the finiteness of the vertex set is not used in the proof. In fact, the relevant chain homotopies are defined on basis tuples and extended linearly, and hence are well-defined for arbitrary graphs. See also Bottinelli--Kaiser \cite[Proposition~3.15, Theorem~3.13]{BK}, where the corresponding contracting homotopy and excision argument are carried out for arbitrary metric spaces.

For a projecting decomposition $(X; G, H)$, we denote inclusions $G\cap H \hookrightarrow G$ and $G \hookrightarrow X$ by $\iota^G$ and $\iota_G$, respectively, and also analogously  for $H$. We also note that we use the convention for the projecting map as $\pi : G\to G\cap H$, while Hepworth--Willerton \cite{HW} use $\pi : H \to G\cap H$ in the original notation. 

\begin{prp}\cite[Theorem~6.5, 6.6]{HW}\label{HWMV}
Let $X$ be a graph, and let $(X; G, H)$ be a projecting decomposition. Then, for all $x \in X$ and all $n, \ell\geq0$, we have the exact sequence
\[
0 \to {\rm MH}^x_{n, \ell}(G\cap H) \xrightarrow{(\iota^G_{\ast}, -\iota^H_{\ast})}  {\rm MH}^x_{n, \ell}(G) \oplus  {\rm MH}^x_{n, \ell}(H) \xrightarrow{\iota_{G\ast}+\iota_{H\ast}}  {\rm MH}^x_{n, \ell}(X) \to 0,
\]
where we define ${\rm MH}^x_{n, \ell}(U) = 0$ if $x \not\in U$ for a subgraph $U\subset X$. Moreover, for $x\in G\cap H$, this exact sequence splits by the retraction
\[
{\rm MH}^x_{n, \ell}(G) \oplus  {\rm MH}^x_{n, \ell}(H)  \xrightarrow{\pi_\ast\oplus 0} {\rm MH}^x_{n, \ell}(G\cap H).
\]
\end{prp}
\begin{proof}
By Remark \ref{Delta}, we may work with the normalized magnitude chain complex ${\rm MC}^{\Delta}_{\ast, \ast}$, which we denote by ${\rm MC}_{\ast, \ast}$ in the following. It is immediate to verify that  the following is an exact sequence of chain complexes
\[
0 \to {\rm MC}^x_{\ast, \ell}(G\cap H) \xrightarrow{(\iota^G_{\ast}, -\iota^H_{\ast})}  {\rm MC}^x_{\ast, \ell}(G) \oplus  {\rm MC}^x_{\ast, \ell}(H) \xrightarrow{\iota_{G\ast}+\iota_{H\ast}}  {\rm MC}^x_{\ast, \ell}(G, H) \to 0,
\]
where we define
\[
{\rm MC}^x_{\ast, \ell}(G, H) = {\rm MC}^x_{\ast, \ell}(G)+{\rm MC}^x_{\ast, \ell}(H) \subset {\rm MC}^x_{\ast, \ell}(X).
\]
Then it is obvious that the projecting map $\pi : G \to G\cap H$ defines a retraction when $x \in G\cap H$.

To show that the inclusion ${\rm MC}_{\ast, \ell}(G, H) \to {\rm MC}_{\ast, \ell}(X)$ is a quasi-isomorphism,  Hepworth--Willerton \cite{HW} constructed a filtration $\mathcal{F}_\bullet$ on ${\rm MC}_{\ast, \ell}(X)$ with 
\begin{align*}
\mathcal{F}_0{\rm MC}_{\ast, \ell}(X) &= {\rm MC}_{\ast, \ell}(G, H), \\
\mathcal{F}_\ell{\rm MC}_{\ast, \ell}(X) &= {\rm MC}_{\ast, \ell}(X),
\end{align*}
and $\mathcal{F}_i{\rm MC}_{\ast, \ell}(X)$ is spanned by tuples $(x_0, \dots, x_n)$ such that  $x_0, \dots, x_{n-i}$ either all lie in $G$ or all lie in $H$. Then they showed each $\mathcal{F}_i/\mathcal{F}_{i-1}$ is acyclic, which implies the desired quasi-isomorphism. From the construction, we have 
\[
\mathcal{F}_i {\rm MC}_{\ast, \ell}(X) \cong \oplus_x \mathcal{F}_i{\rm MC}^x_{\ast, \ell}(X),
\]
for all $i$, where
\[
\mathcal{F}_i {\rm MC}^x_{\ast, \ell}(X) = \mathcal{F}_i {\rm MC}_{\ast, \ell}(X) \cap {\rm MC}^x_{\ast, \ell}(X).
\]
Since $\mathcal{F}_i/\mathcal{F}_{i-1} {\rm MC}_{\ast, \ell}(X)$ is acyclic, $\mathcal{F}_i/\mathcal{F}_{i-1} {\rm MC}^x_{\ast, \ell}(X)$ is also acyclic for all $x \in X$. Hence we obtain that the inclusion ${\rm MC}^x_{\ast, \ell}(G, H) \to {\rm MC}^x_{\ast, \ell}(X)$ is a quasi-isomorphism for all $x \in X$. This completes the proof.
\end{proof}

\begin{prp}\label{BMMV}
Let $X$ be a locally finite graph, and let $(X; G, H)$ be a projecting decomposition. Then, we have the exact sequence
\[
0 \to {\rm MH}^{\rm BM}_{n, \ell}(G\cap H) \xrightarrow{(\iota^G_{\ast}, -\iota^H_{\ast})}  {\rm MH}^{\rm BM}_{n, \ell}(G) \oplus  {\rm MH}^{\rm BM}_{n, \ell}(H) \xrightarrow{\iota_{G\ast}+\iota_{H\ast}}  {\rm MH}^{\rm BM}_{n, \ell}(X) \to 0.
\]
Also, when $X$ has bounded geometry, we have the exact sequence
\[
0 \to {\rm MH}^{\rm uf}_{n, \ell}(G\cap H) \xrightarrow{(\iota^G_{\ast}, -\iota^H_{\ast})}  {\rm MH}^{\rm uf}_{n, \ell}(G) \oplus  {\rm MH}^{\rm uf}_{n, \ell}(H) \xrightarrow{\iota_{G\ast}+\iota_{H\ast}}  {\rm MH}^{\rm uf}_{n, \ell}(X) \to 0.
\]
Moreover, both sequences split by taking the products (or restricted products) of the componentwise retractions with respect to the decomposition in Proposition \ref{prod} induced by the projecting map $\pi : G\to G\cap H$.
\end{prp}
\begin{proof}

By taking product of the exact sequence in Proposition \ref{HWMV}, Lemma \ref{ab4} and Proposition \ref{prod} imply the exact sequence for Borel--Moore magnitude homology. When $X$ has bounded geometry, for fixed $n$ and $\ell$, there is a uniform bound on the ranks of all chain groups involved. With respect to the tuple bases, the boundary, inclusion, and retraction maps are represented by matrices whose entries belong to a fixed finite set. Hence, up to permutations of the tuple bases, only finitely many local chain-complex diagrams occur. Since such permutations preserve the $\ell^\infty$-norms, the induced morphisms on homology have only finitely many norm-preserving isomorphism types. Each induced homology morphism is strict by Lemma \ref{freeabel} (after identifying the relevant subquotients with quotients of finite rank free abelian groups), and the finite number of norm-preserving diagram types makes the strictness uniform by Lemma \ref{fag}. Then Lemma \ref{normedab4}  implies the exact sequence for uniformly finite magnitude homology. The same finite-type argument shows that the family of retractions is uniformly Lipschitz, so their restricted product gives the splitting in the uniformly finite case. This completes the proof.

\end{proof}

\subsection{Mayer--Vietoris sequence for path homology}\label{sec6.2}

\begin{prp}\label{PHMV}
Let $X$ be a locally finite graph, and let $(X; G, H)$ be a projecting decomposition. Then we have a long  exact sequence
\[
 \cdots \to {\rm PH}^{\rm BM}_{n+1}(X) \to {\rm PH}^{\rm BM}_{n}(G\cap H) \xrightarrow{(\iota^G_{\ast}, -\iota^H_{\ast})}  {\rm PH}^{\rm BM}_{n}(G) \oplus  {\rm PH}^{\rm BM}_{n}(H) \xrightarrow{\iota_{G\ast}+\iota_{H\ast}}  {\rm PH}^{\rm BM}_{n}(X) \to \cdots.
\]
Also, when $X$ has bounded geometry, we obtain the same long exact sequence for ${\rm PH}^{\rm uf}_{\ast}$.
\end{prp}
\begin{proof}
Note that the inclusions induce homomorphisms of chain complexes
\[
C_{\ast}(G\cap H) \xrightarrow{(\iota^G_{\ast}, -\iota^H_{\ast})}  C_{\ast}(G) \oplus  C_{\ast}(H) \xrightarrow{\iota_{G\ast}+\iota_{H\ast}}  C_{\ast}(X),
\]
which preserve the filtrations. Hence these inclusions induce morphisms between spectral sequences. By Proposition \ref{BMMV}, we have the following short exact sequence
\[
0 \to {\rm MH}^{\rm BM}_{n, \ell}(G\cap H) \xrightarrow{(\iota^G_{\ast}, -\iota^H_{\ast})}  {\rm MH}^{\rm BM}_{n, \ell}(G) \oplus  {\rm MH}^{\rm BM}_{n, \ell}(H) \xrightarrow{\iota_{G\ast}+\iota_{H\ast}}  {\rm MH}^{\rm BM}_{n, \ell}(X) \to 0,
\]
where $(\iota^G_{\ast}, -\iota^H_{\ast})$ and $\iota_{G\ast}+\iota_{H\ast}$ commute with the differential $d^1$ of the spectral sequence since they are part of a morphism of the spectral sequence. Hence, restricting to the $q=0$ rows, we obtain a short exact sequence of chain complexes with respect to the $d^1$-differentials. Taking homology and using Proposition \ref{pathhomology}, we obtain the desired long exact sequence for Borel--Moore path homology.  It is similar for ${\rm PH}^{\rm uf}_{\ast}$. This completes the proof.
\end{proof}

\begin{rem}
Unlike the case of magnitude homology, the long exact sequence of Proposition \ref{PHMV} need not decompose into short exact sequences. Let us consider the Cayley graph $X$ of $(\Z, \{\pm 1\})$, and its subgraphs $G = \Z_{\leq0}$ and $H = \Z_{\geq0}$. Then it is easy to verify that $(X; G, H)$ is a projecting decomposition.  Note that the function on a one point defines a non-trivial homology class of $ {\rm PH}^{\rm BM}_{0}(G\cap H) \cong {\rm PH}^{\rm uf}_{0}(G\cap H)  $. However, it is a  boundary of a 1-chain both in $G$ and $H$. For example, the boundary of the constant function $1$ on all $(+1)$-edges of $G$ is the point in $G\cap H$, and similarly in $H$. It implies that the induced homomorphism $(\iota^G_{\ast}, -\iota^H_{\ast})$ is not injective at $0$-degree.
\end{rem}

\subsection{Mayer--Vietoris sequence for $E^r_{p, q}$}\label{sec6.3}

\begin{prp}\label{EMV}
Let $X$ be a locally finite graph, and let $(X; G, H)$ be a projecting decomposition. We suppose that each fiber of the projecting map $\pi$ is finite. Then we have the exact sequence
\[
0 \to E^{{\rm BM}, r}_{p, q}(G\cap H) \xrightarrow{(\iota^G_{\ast}, -\iota^H_{\ast})}  E^{{\rm BM}, r}_{p, q}(G) \oplus  E^{{\rm BM}, r}_{p, q}(H) \xrightarrow{\iota_{G\ast}+\iota_{H\ast}}  E^{{\rm BM}, r}_{p, q}(X) \to 0,
\]
for all $r\geq1$ and all $p, q$. Moreover, this exact sequence splits. Also, when $X$ has bounded geometry, we obtain the same split short exact sequence for $E^{{\rm uf}, r}_{p, q}$ provided that the  projecting map $\pi$ has uniformly finite fibers. Note that $G$ and $G\cap H$ also have bounded geometry since they are convex subgraphs of $X$ from the definition of the projecting decomposition and Lemma \ref{Gconv}.
\end{prp}
\begin{proof}
Note that the inclusions induce homomorphisms of chain complexes
\[
C_{\ast}(G\cap H) \xrightarrow{(\iota^G_{\ast}, -\iota^H_{\ast})}  C_{\ast}(G) \oplus  C_{\ast}(H) \xrightarrow{\iota_{G\ast}+\iota_{H\ast}}  C_{\ast}(X),
\]
which preserve the filtrations. Also, from the finiteness assumption on the fiber of $\pi$, it induces a chain map $\pi_\ast : C_{\ast}(G) \to C_{\ast}(G\cap H)$ by Proposition \ref{functorial}, which preserves the filtration since it is 1-Lipschitz.  This defines a retraction $\pi_\ast\oplus0 : C_{\ast}(G) \oplus  C_{\ast}(H) \to C_{\ast}(G\cap H)$ of $(\iota^G_{\ast}, -\iota^H_{\ast})$. Hence $\pi$ induces morphisms between spectral sequences. By Proposition \ref{BMMV}, the sequence on the $E^1$-page is short exact. The morphism induced by $\pi_\ast\oplus0$ is a retraction of its first map, and it commutes with $d^1$. Thus the sequence is a split short exact sequence of chain complexes with respect to $d^1$. Taking homology gives a split short exact sequence on the $E^2$-page. Inductively, if the sequence is split exact on the $E^r$-page, the retraction induced by $\pi_\ast$ commutes with $d^r$, and taking homology yields a split short exact sequence on the $E^{r+1}$-page. It is similar for $E^{{\rm uf}, r}_{p, q}$. This completes the proof.
\end{proof}

\subsection{Mayer--Vietoris sequence for coarse and uniformly finite homology}\label{sec6.4}
From the magnitude--path spectral sequences, we can obtain Mayer--Vietoris sequences for $H_\ast(X)$ and $H^{\rm uf}_\ast(X)$ under the same condition as Proposition \ref{EMV} except for that we do not assume any finiteness on the fiber of the map $\pi$.

\begin{prp}\label{MVCO}
Let $X$ be a locally finite graph, and let $(X; G, K)$ be a projecting decomposition. Then we have the following long exact sequences
\[
\cdots \to H_p(G\cap K) \xrightarrow{(\iota^G_{\ast}, -\iota^K_{\ast})}  H_{p}(G) \oplus  H_{p}(K) \xrightarrow{\iota_{G\ast}+\iota_{K\ast}}  H_{p}(X) \to H_{p-1}(G\cap K) \to \cdots,
\]
and
\[
\cdots \ \to H^{\rm uf}_p(G\cap K) \xrightarrow{(\iota^G_{\ast}, -\iota^K_{\ast})}  H^{\rm uf}_{p}(G) \oplus  H^{\rm uf}_{p}(K) \xrightarrow{\iota_{G\ast}+\iota_{K\ast}}  H^{\rm uf}_{p}(X) \to H^{\rm uf}_{p-1}(G\cap K) \to \cdots,
\]
where the latter is defined when $X$ has bounded geometry.
\end{prp}
\begin{proof}
We give a proof for uniformly finite homology, and the proof for coarse homology is similar. Let $M_\ast = C^{\rm uf}_\ast(G) \oplus C^{\rm uf}_\ast(K) \oplus C^{\rm uf}_{\ast-1}(G\cap K) $ be the mapping cone of the chain map
\[
C^{\rm uf}_\ast (G\cap K) \xrightarrow{(\iota^G_{\ast}, -\iota^K_{\ast})}  C^{\rm uf}_\ast(G) \oplus C^{\rm uf}_\ast(K),
\]
where the differential on $M_\ast$ is defined by $\begin{pmatrix}x\\ y\\ z\end{pmatrix} \mapsto \begin{pmatrix}\del^G_\ast x + \iota^G_\ast z \\ \del^K_\ast y - \iota^K_\ast z\\ -\del^{G\cap K}_\ast z\end{pmatrix} $ for $x \in C^{\rm uf}_\ast(G), y\in C^{\rm uf}_\ast(K)$ and $z \in C^{\rm uf}_{\ast-1}(G\cap K) $. We define a filtration $F_\bullet$ on $M_\ast$ by $F_\ell M_p = F_\ell C^{\rm uf}_p(G) \oplus F_\ell C^{\rm uf}_p(K) \oplus F_\ell C^{\rm uf}_{p-1}(G\cap K) $. Then the inclusion $C^{\rm uf}_\ast(G) \oplus C^{\rm uf}_\ast(K) \to M_\ast$ and the projection $M_\ast \to C^{\rm uf}_{\ast-1} (G\cap K)$ preserve the filtrations. Hence the short exact sequence
\begin{align}\label{long}
0 \to C^{\rm uf}_\ast(G) \oplus C^{\rm uf}_\ast(K) \to M_\ast \to C^{\rm uf}_{\ast-1} (G\cap K) \to 0
\end{align}
induces a short exact sequence
\[
0 \to F_p/F_{p-1}(C^{\rm uf}_\ast(G) \oplus C^{\rm uf}_\ast(K)) \to F_p/F_{p-1}M_\ast \to F_p/F_{p-1}C^{\rm uf}_{\ast-1} (G\cap K) \to 0,
\]
by the nine lemma, and it induces a long exact sequence
\[
\cdots \to E^{{\rm uf}, 1}_{p, q}(G) \oplus E^{{\rm uf}, 1}_{p, q}(K) \to E^{{\rm uf}, 1}_{p, q}(M) \to E^{{\rm uf}, 1}_{p, q-1}(G\cap K)\xrightarrow{\delta_\ast} E^{{\rm uf}, 1}_{p, q-1}(G) \oplus E^{{\rm uf}, 1}_{p, q-1}(K) \to \cdots , 
\]
for all $p, q$. The standard argument shows that the connecting homomorphism $\delta_\ast$ above coincides with the homomorphism induced from the chain map $(\iota^G_{\ast}, -\iota^K_{\ast})$, which is injective by Proposition \ref{BMMV}. Hence we obtain a short exact sequence
\[
0\to E^{{\rm uf}, 1}_{p, q}(G\cap K) \to E^{{\rm uf}, 1}_{p, q}(G) \oplus E^{{\rm uf}, 1}_{p, q}(K) \to E^{{\rm uf}, 1}_{p, q}(M)  \to 0.
\]
Moreover, the chain map $\iota_{G\ast} \oplus \iota_{K\ast}\oplus 0 : M_\ast \to C^{\rm uf}_\ast(X)$ induces 
a commutative ladder
\begin{equation*}
\xymatrix{
0\ar[r]  &  E^{{\rm uf}, 1}_{p, q}(G\cap K) \ar[r]\ar@{=}[d] & E^{{\rm uf}, 1}_{p, q}(G) \oplus E^{{\rm uf}, 1}_{p, q}(K) \ar[r]\ar[r]\ar@{=}[d] & E^{{\rm uf}, 1}_{p, q}(M) \ar[r]\ar[d] & 0 \\
0\ar[r] &  E^{{\rm uf}, 1}_{p, q}(G\cap K) \ar[r] & E^{{\rm uf}, 1}_{p, q}(G) \oplus E^{{\rm uf}, 1}_{p, q}(K) \ar[r] & E^{{\rm uf}, 1}_{p, q}(X) \ar[r]& 0,
}
\end{equation*}
hence this chain map induces an isomorphism on the $E^1$-page by the five lemma. Thus, since the filtrations on $M_\ast$ and $C_\ast^{\rm uf}(X)$ are bounded below and exhaustive, Proposition \ref{comparison} implies the isomorphism $H_\ast(M) \cong H^{\rm uf}_\ast(X)$, and hence we obtain the desired long exact sequence from  the short exact sequence (\ref{long}). This completes the proof.
\end{proof}

\begin{rem}\label{mitchener}
The Mayer--Vietoris sequence for coarse homology in Proposition \ref{MVCO} is closely related to the general coarse Mayer--Vietoris principle. In the axiomatic formulation of coarse homology theories, a coarsely excisive decomposition gives rise to a Mayer--Vietoris long exact sequence; see, for example, \cite[Definition~3.2 in addendum]{Mit}. In fact, the excisive condition $\forall R>0, \exists S>0, N_R(G)\cap N_R(H)\subset N_S(G\cap H)$ is always satisfied when $G, H$ are subgraphs of a graph $X$ and $G\cup H = X$. The argument above gives a direct chain-level proof in our setting via the magnitude--path spectral sequence. We do not know whether this uniformly finite version has appeared explicitly in the literature.
\end{rem}

\begin{rem}
It is natural to ask whether {\it the K\"{u}nneth formula} holds for coarse or uniformly finite homology. In the setting of the ordinary magnitude--path spectral sequence, Hepworth--Roff proved  K\"{u}nneth theorems for its pages with respect to the box product \cite[Theorem 5.6]{HRbigraded}; see also \cite[Theorem 5.2]{HR} for the corresponding result for reachability homology. Their proof is based on a filtered Eilenberg--Zilber-type chain homotopy equivalence ${\rm RC}_\ast(G)\otimes {\rm RC}_\ast(H)\simeq {\rm RC}_\ast(G\square H)$\cite[Theorem~4.7]{HRbigraded}. 

It is therefore natural to ask whether an analogous argument can be developed for the Borel--Moore and uniformly finite magnitude--path spectral sequences considered here. The infinite-chain versions considered here are more subtle, however. Already in degree zero, neither $C_\ast(X\times Y)$ nor $C_\ast^{\rm uf}(X\times Y)$ can in general be identified with the corresponding algebraic tensor product. Thus the usual Eilenberg--Zilber and algebraic K\"{u}nneth arguments do not immediately apply. Product structures in uniformly finite homology have been studied by, for example,  Diana \cite[Chapter~3]{Diana} and Diana--Nowak \cite{DN}.  It would be interesting to investigate whether the magnitude--path spectral sequence provides a useful approach to this problem.
\end{rem}

\end{document}